%Format: plain
\input amstex
\input amsppt.sty

%
%                On real structures of complex surfaces
%                      and applications
%
%                       Viatcheslav Kharlamov, Viktor Kulikov
%

\pagewidth{32.5pc} \pageheight{
%43
48pc} \topskip=\normalbaselineskip

\def\stydate{November 27, 1997}
\immediate\write16{This is DEGT.DEF by A.Degtyarev <\stydate>}
{\edef\temp{\the\everyjob\immediate\write16{DEGT.DEF: <\stydate>}}
\global\everyjob\expandafter{\temp}}

\chardef\tempcat\catcode`\@\catcode`\@=11

%       Redefine \ge and \le:
\let\ge\geqslant
\let\le\leqslant
%       Frequently used \Bbb letters:
\def\C{{\Bbb C}}
\def\R{{\Bbb R}}
\def\Z{{\Bbb Z}}
\def\Q{{\Bbb Q}}
\def\N{{\Bbb N}}

%       Notation for projective spaces (usage: \Rp3):
\def\Cp#1{\C{\operator@font p}^{#1}}
\def\Rp#1{\R{\operator@font p}^{#1}}
%       Plenty of operators not covered by PLAIN.TEX:
\def\Alb{\qopname@{Alb}}
\def\Hom{\qopname@{Hom}}
\def\Ext{\qopname@{Ext}}
\def\Tors{\qopname@{Tors}}

\def\Im{\qopname@{Im}}          %   I use these
\def\Re{\qopname@{Re}}     %
\def\Ker{\qopname@{Ker}}
\def\Coker{\qopname@{Coker}}
\def\Int{\qopname@{Int}}
\def\Cl{\qopname@{Cl}}
\def\Fr{\qopname@{Fr}}
\def\Fix{\qopname@{Fix}}
\def\tr{\qopname@{tr}}
\def\inj{\qopname@{in}}
\def\id{\qopname@{id}}
\def\pr{\qopname@{pr}}
\def\rel{\qopname@{rel}}
\def\pt{{\operator@font{pt}}}
\def\const{{\operator@font{const}}}
\def\codim{\qopname@{codim}}
\def\cdim{\qopname@{dim_{\C}}}
\def\rdim{\qopname@{dim_{\R}}}
\def\conj{\qopname@{conj}}
\def\rank{\qopname@{rk}}
\def\sign{\qopname@{sign}}
\def\gcd{\qopname@{g.c.d.}}

\def\set<#1|#2>{\bigl\{#1\bigm|#2\bigr\}}

%       This is to be used after \\ in \CD, \align..., and \matrix

%       Some proofreading tools:
\def\preprint#1{\hrule height0pt depth0pt\kern-24pt%
  \hbox to\hsize{#1}\kern24pt}
\def\today{\ifcase\month\or January\or February\or March\or
  April\or May\or June\or July\or August\or September\or October\or
  November\or December\fi \space\number\day, \number\year}

\def\n@te#1#2{\leavevmode\vadjust{%
 {\setbox\z@\hbox to\z@{\strut\eightpoint#1}%
  \setbox\z@\hbox{\raise\dp\strutbox\box\z@}\ht\z@=\z@\dp\z@=\z@%
  #2\box\z@}}}
\def\leftnote#1{\n@te{\hss#1\quad}{}}
\def\rightnote#1{\n@te{\quad\kern-\leftskip#1\hss}{\moveright\hsize}}
\def\?{\FN@\qumark}
\def\qumark{\ifx\next"\DN@"##1"{\leftnote{\rm##1}}\else
 \DN@{\leftnote{\rm??}}\fi{\rm??}\next@}

\def\centerpage{\dimen@=6.5truein \advance\dimen@-\hsize\hoffset.5\dimen@}
\ifnum\mag>1000 \centerpage\fi

%       To get rid of AmSTeX's logo:
\def\nologo{\let\logo@\relax}

%  This is done for compatibility with AmS-LaTeX
\expandafter\ifx\csname eat@\endcsname\relax\def\eat@#1{}\fi
\expandafter\ifx\csname operator@font\endcsname\relax
 \def\operator@font{\roman}\fi
\expandafter\ifx\csname eightpoint\endcsname\relax
 \let\eightpoint\small\fi

\catcode`\@\tempcat\let\tempcat\undefined

\hoffset1in                 %%%%%%% for proofreading

\let\mnote\leftnote
\def\longnote#1#2#3{\mnote{#2 \ #1}\bgroup\eightpoint
 \plainfootnote{}{\?"#1" \ #3}\egroup}
 \hoffset.5in

\let\mnote\leftnote

\let\+\sqcup
\def\er{E_{\R}\futurelet\nexts\getsup}
\def\xr{X_{\R}\futurelet\nexts\getsup}

\def\c{\futurelet\nexts\conjs}
\def\index#1{^{\botsmash{(#1)}}}
\def\conjs{\ifx\nexts1\def\nexts##1{t\index1}\else
 \ifx\nexts2\def\nexts##1{t\index2}\else
 \def\nexts{\mathop{\roman{conj}}}\fi\fi\nexts}

\def\getsup{\ifx\nexts1\def\nexts##1{\index1}\else
 \ifx\nexts2\def\nexts##1{\index2}\else\let\nexts\relax\fi\fi\nexts}

\def\rel{{\operatorname{rel}}}

\def\li{l}

\def\<#1>{\langle#1\rangle}

%%%%%%%%%%%%%%%  Definitions for Figures 1, 2. %%%%%%%%%%%%%%%
\newdimen\step \step10pt
\newdimen\thin \thin.2pt
\newdimen\thick \thick.4pt
\newdimen\gap \gap1pt
\newcount\loopcount

\def\grid#1#2{\vbox{\offinterlineskip \parindent 0pt \hsize#1\step
 \advance\hsize by\thick \advance\hsize by\gap
 \dimen0 \hsize \kern\gap \advance\hsize-\thin
 \def\vert##1{\vbox to0pt{\hsize##1 %
  \vss\vrule height\dimen0 width##1}\kern-\thin}%
 \def\zbox##1{\vbox to0pt{\hsize0pt \vss\hbox{\raise##1 \box0 }}}
 \def\labelbox##1{\setbox0 \hbox to0pt{$\scriptstyle##1$}%
  \dimen1=-\ht0 \ht0=0pt \dp0=0pt}
 \def\xlabel##1{\line{\kern\hsize\labelbox{\,\,##1\hss}\zbox{1pt}\hss}}
 \def\ylabel##1{\line{\labelbox{\hss##1\,}\advance\dimen1 by\dimen0 %
  \zbox{\dimen1}\hss}}
 \def\pointbox##1(##2,##3){\dimen1\step \multiply\dimen1 by##2 %
  \dimen2\step \multiply\dimen2 by##3
  \setbox0\vbox to0pt{\hsize 0pt\vss\kern-2\thick\kern2\thin\line
   {\hss\kern2\thick\kern2\thin$##1$\hss}\vss}%
  \moveright\dimen1\hbox{\zbox{\dimen2}}}
 \def\black{\pointbox\bullet}
 \def\white{\pointbox\circ}
 \loopcount#1 \loop\hrule height\thin width\hsize \kern\step
  \kern-\thin \advance\loopcount-1 \ifnum\loopcount>0 \repeat
 \hrule height\thick \line{\vert\thick
 \loopcount#1 \loop \kern\step \vert\thin
  \advance\loopcount-1 \ifnum\loopcount>0 \repeat \hfil}
 \line{\vbox{#2}}}}

\def\M{\Cal M}

\def\letitem#1{\item"\setbox0\hbox{(x)}\hbox to\wd0{(#1)\hss}"}

\def\Aut{\operatorname{Aut}}

\def\dif{\operatorname{Dif}}
\def\iso{\operatorname{Iso}}
\def\Kl{\operatorname{Kl}}
\def\Def{\operatorname{Def}}
\def\kl{\operatorname{kl}}

\catcode`\@=11
\def\smallparn#1{\smash{\vcenter{\hbox{$\m@th\scriptscriptstyle#1$}}}}
\def\smallscript#1\\{\mathopen{\smallparn(}#1\mathclose{\smallparn)}}

\catcode`\@\active

\def\V(#1){V_{#1}}

%%%%%%%%%%%%%%%%%%%%   Two column display %%%%%%%%%%%%%%%%%%%%%%%%%%%
\newdimen\twocoljot \twocoljot\jot
\newdimen\firstindent \firstindent2\parindent
\newdimen\secondindent \secondindent\parindent
\def\twocol{\par\bgroup\interlinepenalty200%\interdisplaylinepenalty
 \openup\twocoljot\halign\bgroup\hbox to.5\hsize
  {\kern\firstindent$##$\hfil}&\kern\secondindent$##$\hfil\cr}
\def\endtwocol{\crcr\egroup\egroup\par\smallpagebreak}
\def\label(#1){\llap{\rm(#1)\space}}

\let\mnote\leftnote

\topmatter

\title    Deformation inequivalent complex conjugated
complex structures
and applications
\endtitle
\rightheadtext{
Deformation inequivalent complex structures} %   Short title

\author
Viatcheslav Kharlamov and Viktor Kulikov$^*$
\endauthor

\address
Institut de Recherche Math\'ematique Avanc\'ee\newline\indent
Universit\'e Louis Pasteur et CNRS\newline\indent
7 rue Ren\'e-Descartes
\newline\indent 67084, Strasbourg, France
\endaddress

\email
kharlam\,\@\,math.u-strasbg.fr
\endemail

\address
Steklov Mathematical Institute\newline\indent
\newline\indent Moscow, Russia
\endaddress

\email
kulikov\,\@\,mi.ras.ru
\endemail



\thanks \flushpar $^*$
Partially supported by INTAS-00-0259,
NWO-RFBR-047-008-005, and
RFBR  99-01-01133.
\endthanks

\abstract
We start from a short summary of our principal result from \cite{KK}:
an example of a complex algebraic surface which is
not deformation equivalent to its complex conjugate and
which, moreover, has no homeomorphisms reversing the canonical class.
Then,
we generalize this
result to higher dimensions and construct several series of higher
dimensional compact complex manifolds having no homeomorphisms reversing
the canonical class.
After that,
we resume and broaden 
the 
applications given in \cite{KK} and
\cite{KK2}, in particular, as a new application, we propose examples of
(deformation) non equivalent symplectic structures with opposite
canonical classes. 
\endabstract

\endtopmatter

\document

\head
1.\enspace
Introduction
\endhead

Many of achievements in real algebraic geometry appeared as applications
of complex geometry.
By contrary, the results of the present paper, which has grown
from a solution of some questions from the real algebraic geometry,
can be considered as applications in the opposite direction.
They may be 
of an interest for symplectic geometry too.

To state the principal questions we need to fix some definitions.
We choose the language of complex analytic geometry since in 
this 
setting the results sometimes look stronger than if
we restricted ourselves to algebraic or K\"ahler
varieties, though the K\"ahler hypothesis could simplify several proofs and
could allow one to extend the set of examples.

Thus, we define a {\it real
structure \/} on a complex manifold $X$ as an antiholomorphic
involution $c:X\to X$.
By a {\it
deformation equivalence \/} we mean the
equivalence generated by local
deformations of complex manifolds.
Since we are interested only in compact
manifolds, this, commonly used,
local deformation equivalence relation
can be defined in the following way:
given a proper holomorphic submersion
$f:W\to B_1,\, B_1=\{\vert z \vert < 1\}\subset \C$, the complex
varieties {\bf isomorphic}
to $X_t=f^{-1}(t), t\in B_1$,
are said to be {\it deformation equivalent}\rom.

Note that anti-holomorphic automorphisms ({\it anti-automorphisms}
for short) of a complex manifold $X=(M,J)$, where $M$ is the underlying
smooth manifold and $J$ is a complex structure
on it, can be
regarded 
as isomorphisms between $X$ and its complex
conjugated twin, $\bar X=(M,-J)$. All the anti-automorphisms together
with the automorphisms form a group which we denote by $\Kl$ and call
{\it Klein group}. Clearly, the Klein group is either a $\Z/2$-extension
of the automorphism group
$\Aut$, $1\to\Aut\to\Kl
@>\kl>>
\Z/2\to 1,$ or
coincides with it, $\Kl=\Aut$.

Note also that the underlying smooth manifolds of deformation equivalent
varieties are necessary diffeomorphic. Moreover, the diffeomorphisms
coming from a deformation preserve the complex orientation and the
canonical class. Certainly, in the same time,
on a given differentiable manifold here may exist 
isomorphic complex structures with different canonical classes or
even with different complex orientations.

Below are two closely related questions from real
algebraic geometry, which initiated our study.
{\it
Can any complex compact manifold be deformed to a
manifold which can be equipped with a real structure or at least to a
manifold with an anti-automorphism} (not necessarily of order $2$)?

It is clear that the response is in
the affirmative as
long as the (complex) dimension of the manifold is $\le 1$, i.e., for
points and Riemann surfaces: obviously, any Riemann surface can be
deformed to a real one. To our knowledge, starting from $\dim_\C=2$ the
both questions remained open. We have proved the following theorem.

\proclaim{Theorem 1.1}
In any dimension $\ge 2$ there
exists a compact complex manifold $X$ which
can not be deformed
to $\bar X$\rom; in particular, $X$ cannot be deformed to a real manifold
or to a manifold with an anti-automorphism.
\endproclaim

Explicit examples are given below in Sections 2 and 3.
Many of them have a stronger property (Corollary 3.13):

{\it In any dimension $\ge 2$ there exists a projective
complex manifold $X$ of general type
which has no homeomorphisms $h:X\to X$ such that
$h^*c_1(X)=-c_1(X)$.}

It is worth noticing that in dimension $2$
all our examples are rigid, i.e., they have no
local nontrivial deformations at all.
Moreover, they are strongly rigid, i.e., up to isomorphism
and conjugation, they have only one complex structure.
Their strong rigidity is not used in the proof of Theorem 1.1
in its dimension $2$ part, but we use
it in higher dimensions.

The examples mentioned above
provide the promised applications, which are: new counter-examples
to the $\dif=\Def$ problem for
complex surfaces (and higher dimensional varieties);
first (to our knowledge) 
counter-examples to the
ambient $\dif\Rightarrow\Def$ problem for plane cuspidal curves; and
first (to our knowledge) examples of symplectic structures $\omega$ not
equivalent, up to diffeomorphisms and deformations, to their reverse 
$-\omega$.

Deformation equivalent complex manifolds are orientedly ($\Cal
C^\infty$-) diffeomorphic and it is the converse which is called the
$\dif=\Def$ {\it problem} (the, probably, first mentioning of this
problem is found in \cite{FM}). Obviously, the response is in 
the affirmative
in $\dim\le 1$. For complex
surfaces the question remained open for quit a
while. The first counter-examples were found by Manetti only a couple of
years ago, see \cite{Ma}. The surfaces $X$ and $\bar X$ from our Theorem
1.1 are different and, by our opinion, more simple counter-examples. At
least, they have the following, useful for our further applications,
properties: our inequivalent complex structures are opposite to each
other ($J_0=-J_1$; in other words, the complex structures are complex
conjugated), so that their canonical classes are opposite ($c_1(J_0)=-
c_1(J_1)$), and, as in Corollary 3.13, there is no homeomorphism
transforming $c_1(J_0)$ in $c_1(J_1)$. Moreover, in our examples the
structures $J_0$ and $J_1$ are not equivalent even in the class of
almost-complex structures, see the concluding Remark in Section 4.

The counter-examples to the $\dif=\Def$ problem in $\dim\ge 3$ were known
before. In dimension $3$, at least in the category of K\"ahler manifolds,
one can take products of the Riemann sphere with Dolgachev surfaces, see
\cite{R}; as it follows from the classification of
$6$-dimensional manifolds \cite{Z}, they are diffeomorphic, and their
stability under deformations in the category of K\"ahler manifolds follows
from \cite{Fujiki} (note that, as is proven in \cite{R}, these
counter-examples work as well in the symplectic category). In dimension
$4$ and in higher even dimensions one can take products of Dolgachev
surfaces and argue on divisors of the canonical class which are invariant
under deformation, see \cite{FM}; such products are diffeomorphic due to
$h$-cobordance of Dolgachev surfaces, which follows from their homotopy
equivalence, see \cite{W}. Note that in our $\dim\ge 3$ examples the non
equivalence
of the varieties is not
related to the divisibility of the canonical class and, in each
dimension, among the examples there are varieties of general type.

It is worth noticing also that the
convention to fix the orientation in the
statement of the $\dif=\Def$ problem
is not necessary in the case of
complex dimension $2$. As it follows from Kotschick's results,
see \cite{Ko}, in this dimension
the existence of two complex structures with opposite
orientations (which is an extremely rare phenomenon
in this dimension) implies
the existence of an orientation-reversing diffeomorphism, and, thus,
the set of equivalence classes does not change if the orientation
convention is removed.

Now, let us discuss the second application, which concerns
the {\it ambient $\dif\Rightarrow\Def$ problem for plane curves}\rom.
As in the absolute
case, (equisingular) deformation equivalent pairs $(P^2,C)$,
where $C\subset P^2$ stands for a plane
curve, are diffeomorphic 
whenever $C$ is a cuspidal curve, and the problem consists in the 
converse statement (for non cuspidal curves the diffeotopy should be
replaced by an isotopy; in fact, in all cases one can speak already on
$\dif\Rightarrow\iso$ problem; note that the isotopy coming from an
equisingular deformation certainly always can be made smooth outside the
singular points, and for some types of singular points, as in the case of
ordinary nodes and cusps, the isotopy can be made smooth everywhere). If
the curve is nonsingular or if it has only ordinary nodes, the response
to the $\dif\Rightarrow\Def$ problem is in the affirmative 
(for the nonsingular
case it is trivial, for curves with ordinary nodes it follows from the
Severi assertion on the irreducibility of the corresponding spaces). We
have produced, using our surfaces from theorem 1.1, an infinite sequence
of counter-examples in which the curves are cuspidal. Our curves are
irreducible and in each example the diffeomorphism $P^2\to P^2$
transforming one curve into another is, in fact, the standard complex
conjugation. They are obtained as the visible contours of the surfaces
from theorem 1.1 (see more details in Section 4).

To conclude the introduction, let us mention also an application
to a $\dif=\Def$ {\it problem in symplectic geometry}\rom, which may 
also be of a certain interest. Consider, on a given
oriented smooth $4$-manifold, all the possible symplectic structures
compatible with the orientation. Call two structures {\it equivalent} if 
one can be obtained from the other by a deformation
followed, if necessary, by a diffeomorphism. Our surfaces from theorem
1.1 provide examples where the number of equivalence classes is at least
$2$; namely, in these examples any K\"ahler symplectic form $\omega$ and
its reverse, $-\omega$, are not equivalent to each other; it is due to
the fact that there is no diffeomorphism 
(or even homeomorphism) reversing the
canonical class of $\omega$, see section 4 for details.
To our knowledge, in the previously
known examples of nonequivalent symplectic structures $\omega_1,
\omega_2$ (see \cite{MT}, \cite{LeB}, and \cite{Sm})
the assertion states only that
$\omega_1\ne\pm\omega_2$. %, whatever is the dimension.

It is probably worth to note that in our examples
the nonequivalence can not be detected
by the divisibility properties
of the canonical class or by the Seiberg-Witten invariant: we deal with
surfaces of general type, our complex and symplectic structures are
opposite, as well as their canonical classes, and thus in our examples
the Seiberg-Witten functions coincide.

{\bf Remarks.} In dimension two, in our counter-examples to the
$\dif=\Def$ problem the number of equivalence classes is equal to $2$,
and, moreover, the moduli space is merely a two point set.
%Short after their appearance, F.~Catanese \cite{Ca}
%worked out other 
Additional
dimension two examples
of deformation inequivalent complex
conjugated complex structures
%based on a different idea.
are worked out by F.Catanese in \cite{Ca}.

{\bf Acknowledgements.} The first author is grateful to the participants
of the Gokova conference for their stimulating interest.
We thank also J.~Koll\'ar,
D.~Kotschick, M.~Kreck,
W.-J.~Li, and M.~Paun for
their helpful assistance during the preparation of this expository paper.

\head 2.\enspace Principal examples.
\endhead

Our first, main, example is a
suitably
chosen ramified
Galois $(\Z/5\times\Z/5)$-covering of $\Cp2$
(its construction is inspired by some Hirzebruch's one, cf., \cite{Hirz}).
The ramification locus is the configuration of $9$ lines dual to
the inflection points
of a nonsingular cubic in the dual plane, $(\Cp2)^*$.
What follows does not depend
on a particular choice of the cubic, moreover, even the
configuration of the 9 lines
does not depend, up to projective transformation, from this choice. In proper
chosen homogeneous coordinates it is defined by equation
$(x_1^3-x_2^3)(x_2^3-x_3^3)(x_3^3-x_1^3)=0$
(as it follows, indeed, from the famous elementary geometry  C\'eva
theorem; the equation given
corresponds to the Ferma cubic $x_1^3+x_2^3+x_3^3=0$). In these configuration
there are $12$ triple points and they are the only multiple points.

Birationally, our covering is given by equations
$$
z^5=l_1l_2l_3l_4^3l_5^3l_9,\quad w^5=l_1l_3l_4^3l_6l_7l_8^2l_9,
$$
where the lines $l_1=0,\dots, l_9=0$ are numbered
like in the table below
following the identification of inflection
points with the points of order $3$ in the
Jacobian of the cubic.
In the homogeneous coordinates like above
one can take
$$
\matrix
\li _1=%\{
x_1-x_3, %=0 \},
& \li _2=%\{
x_1-\mu^2 x_3, %=0 \},
& \li _3=%\{
x_1+\mu x_3, %=0 \},
\\
\li _4=%\{
x_2-\mu^2x_3, %=0 \},
& \li _5=%\{
x_2-x_3, %=0 \},
&
\li _6=%\{
x_2+\mu x_3=0, %  \},
\\
\li _7=%\{
x_1+\mu x_2, %=0  \},
& \li _8=%\{
x_1-\mu^2x_2, %=0 \},
&
\li _9=%\{
x_1-x_2. %=0 \}.
\endmatrix
$$
where $\mu =e^{\pi i/3}$.

\centerline{\it %Table of c
Correspondence between the lines and
the points of order 3 in the Jacobian.}
%Format: plain
$$\phantom{aaaa}$$

\centerline{\vbox{%
\def\\{\cr\noalign{\hrule}}
\def\plain#1{\omit\hss#1\hss\srule}
\def\margin{\kern6pt}                 %% Adjust cell margins here
\def\srule{\vrule height9pt depth3pt} %% Adjust row height and depth here
\halign{%
 \vrule\margin\hss#\hss\margin\srule&&    %% First entry centered
 \margin\hss#\margin\vrule\cr           %% Others flushed right
 \noalign{\hrule}
% Table starts here
 $l_1$ & $l_4$ & $l_7$ \\ %% Column headers
 $l_2$ &  $l_5$ &  $l_8$ \\
 $l_3$ & $l_6$  &  $l_9$ \\
% Table ends here
 \crcr}}}
$$
\text{
%Fig.
Table 1}
%\phantom{aaaa}
$$

As any Galois covering,
our one is determined by the corresponding
homomorphism $\phi=(\phi_1,\phi_2) :
H_1(\Cp2\setminus\cup l_r)\to\Z/5\times\Z/5$
(here and further, we do not make difference between
the notation of the lines and their linear form expressions).
The exponents in the equations of the covering
encode the value of $\phi$
on the standard, dual to the lines, generators $\lambda_r\in
H_1(\Cp2\setminus\cup l_r)$:
$$
\matrix \phi(\lambda _1)=(1,1), & \phi(\lambda _2)=
(1,0), & \phi(\lambda _3)=(1,1), \\
\phi(\lambda _4)=(3,3), & \phi(\lambda _5)=
(3,0), & \phi(\lambda _6)=(0,1), \\
\phi(\lambda _7)=(0,1), & \phi(\lambda _8)=
(0,2), & \phi(\lambda _9)=(1,1),%.
\endmatrix
$$

Thus defined $(\Z/5\times\Z/5)$-covering of $\Cp2$, $X\to\Cp2$,
has isolated
singularities, which arise from the $12$ triple points of the
ramification locus. To get its (minimal) resolution it is sufficient
to blow up $\Cp2$ in the triple points and to take
the induced
$(\Z/5\times\Z/5)$-covering of the
blown up plane, $\tilde X\to\Cp2(12)$.

Let denote by $D_{ijk}\subset\tilde X$
the full transform of the
blown up curves $E_{ijk}\subset \Cp2(12)$
and by $C_r\subset\tilde X$
the strict transform
of the lines $l_r\subset \Cp2$.
All of them belong to the ramification locus of
$\tilde X\to\Cp2(12)$ and they all have the ramification index equal to $5$.
For the intermediate strict transforms $L_r\subset\Cp2(12)$
of $l_r$ one has $L_r^2=-3$.
So, elementary pull back calculation shows that
$$
C_r^2=-3, \quad D^2_{ijk}=-1,\quad\text{and}\quad 3K=7\sum C_r+12\sum D_{ijk}.
\tag 1
$$
In addition,
$$
(C_r,K)=9\quad \text{and}\quad (D_{ijk},K)=3.
\tag 2
$$

%The following result is the key point.

\proclaim{Proposition 2.1} The surface $\tilde X$ is
a strongly rigid surface
of general type with ample canonical divisor.
The group $\Kl(\tilde X)$ coincides with the
covering transformations group $G=\Z/5\times\Z/5$.
In particular, there does not exist
neither a real structure nor even an
anti-holomorphic diffeomorphism on $\tilde X$.
\endproclaim

Here and everywhere further,
by the {\it strong} or {\it Mostow-Siu rigidity}
of a compact complex manifold $M$ we mean
the following property:
whatever are a compact complex K\"ahler variety $Y$ and
a
%$C^\infty$-submersion
continious map $p: Y\to M$ with
nonzero $p_*:H_{2m}(Y;\Z)\to H_{2m}(M;\Z)$,
$m=\dim_\C M$, then $p$ is homotopic
to a holomorphic or anti-holomorphic map
$Y\to M.$ According to Siu's results \cite{S}, any non singular compact
quotient of an irreducible bounded Hermitian symmetric domain
of dimension $\geq 2$
has such a property. In addition, as is known, for any such quotient
$M$ one has $H^1(M;\Theta)=0$ ($\Theta$ states for the sheaf
of holomorphic tangent fields),
which implies, in  particular, that
$M$ is {\it locally rigid},
i.e., its complex structure has no non trivial local deformations.

{\it Sketch of the proof of Proposition 2.1} \rom(see \cite{KK} for
calculations and combinatorial details).
A direct calculation shows that
$K_{\tilde X}^2=333$ and $e(\tilde X)=111$
($e$ states for the Euler characteristic), so
$\tilde X$ is a so-called
Miyaoka-Yau surface;
%the fact that
%such surfaces are strongly rigid
their universal cover is  a complex ball,
%and, moreover, unique in their homotopy
%type up to
%holomorphic and anti-holomorphic diffeomorphisms
%is well known,
see, for example,
\cite{BPV},
and, thus, their strong rigidity follows from \cite{S}
(all the complex structures on the underlying smooth manifold are K\"ahler,
as it follows, for example, from $b_1=0$; the
%above uniqueness
strong rigidity applied to the identity map
implies also the absence of complex structures
with opposite orientation).
Due to (1), (2) and the Nakai-Moishezon criterion,
$K$ is ample.

Since the ramification locus contains
the strict transform of at least two
different pairs of lines,
to prove $\Kl=G$ it is sufficient to show, first,
that both the union of strict
transforms $C_r$ of the lines
and the union of the exceptional divisors $D_{ijk}$ are
invariant under Klein transformations of $\tilde X$, and, second,
that there is no Klein transformation of $\Cp2$, different from $\id$,
which can be lifted from $\Cp2$ to $\tilde X$
({\it respecting the covering} for short).

Suppose that $\cup C_r$
is not preserved under some $h\in \Kl$,
i.e., there is $i$ such that $h(C_i)\not\subset C=\cup C_r$. Then,
$(h(C_i), \sum C_r)=a\ge 0$, $(h(C_i), \sum D_{ijk})=b\ge 0$,
and
$3K=7\sum C_r+12 \sum D_{ijk}$ implies
$$
7a+12b=3(h(C_i),K)=27.
$$
To get a contradiction
it remains to note that the latter equation
has no solutions in $a,b\in\Z_+$ (here, counting the intersection
numbers we equip all the curves, including $h(C_i)$,
with their complex orientation induced from $\tilde X$).
The proof of $h(D)\subset D$, $D=\cup D_{ijk}$, is completely similar.

Now, it follows that every $h\in\Kl(\tilde X)$ is lifted
from $\Cp2$ and thus it remains to prove that the only $g\in\Kl(\Cp2)$
respecting the covering is $g=\id$. Since
$g$ respects the covering, it acts on the set of
intermediate $\Z/5$-coverings $Y\to\Cp2$ and their deck transformations.
Namely, such a covering with a marked deck transformation is given
by an epimorphism $\psi : H_1(\Cp2\setminus\cup l_r)\to\Z/5$ and $g$
transforms it into $(-1)^{\kl g}\psi\circ g_*$
(recall that $\kl g$ is $0$ if $g$ is holomoprhic and
$1$ otherwise).
The following is the table of
rows representing all the epimorphisms
$\psi=x\phi_1+y\phi_2 : H_1(\Cp2\setminus\cup l_r)\to\Z/5$
in the base $\lambda_i$:

$$
\matrix
(1,1,1,3,3,0,0,0,1) &(2,2,2,1,1,0,0,0,2) &(3,3,3,4,4,0,0,0,3) &(4,4,4,2,2,0,0,0,4)\\
(1,0,1,3,0,1,1,2,1) &(2,0,2,1,0,2,2,4,2) &(3,0,3,4,0,3,3,1,3) &(4,0,4,2,0,4,4,3,4)\\
(2,1,2,1,3,1,1,2,2) &(4,2,4,2,1,2,2,4,4) &(1,3,1,3,4,3,3,1,1) &(3,4,3,4,2,4,4,3,3)\\
(3,1,3,4,3,2,2,4,3) &(1,2,1,3,1,4,4,3,1) &(4,3,4,2,4,1,1,2,4) &(2,4,2,1,2,3,3,1,2)\\
(4,1,4,2,3,3,3,1,4) &(3,2,3,4,1,1,1,2,3) &(2,3,2,1,4,4,4,3,2) &(1,4,1,3,2,2,2,4,1)\\
(0,1,0,0,3,4,4,3,0) &(0,2,0,0,1,3,3,1,0) &(0,3,0,0,4,2,2,4,0) &(0,4,0,0,2,1,1,2,0).
\endmatrix
$$

The elements in a row represent, in fact, the weights of the
deck transformation on the components of the ramifications locus.
Thus, they should be preserved up to permutation.
In particular, the functions $a\mapsto r_i(a)$
counting the number of coordinates of a row $a$ equal $i$
should be invariant under the action of $g$.

An easy running over, taking into account, in addition,
the incidence relations
between the lines $l_i$, shows that the only possibility for $g$
is to keep invariant each line. Hence, $g=\id$.\qed

Theorem 1.1 in dimension two
is a straightforward consequence of
the second statement of Proposition 2.1.
In fact, 
our proof of this part of Proposition 2.1
uses only the
local rigidity of $\tilde X$, which is an easier result than the
strong rigidity. As to the latter, it implies the following
property which is crucial for the other applications and for the
proof of Theorem 1.1 in higher dimensions.

\proclaim{Lemma 2.2}
The homeotopy group of $X$ \rom(i.e., the group of homotopy classes of
homeomorphisms $X\to X$\rom) coincides with the covering transformations
group $G=\Z/5\Z\times\Z/5\Z$. In particular,
for any homeomorphism
$f:\tilde X\to \tilde X$
one has $f^*[K]=[K], [K]\in H^2(X;\Z)$.
\endproclaim

\demo{Proof} By strong rigidity, any homeomorphism $\tilde X\to \tilde X$
is homotopic to an automorphism or to an anti-automorphism. Due to
Proposition 2.1, there are no anti-automorphisms and the group of
automorphisms coincides with the covering transformations group $G$.
It remains to note that the canonical class is preserved
by any automorphism.
\qed
\enddemo

Another important property of $\tilde X$ is that its irregularity
$q=q(\tilde X)$ is zero, see Appendix.

There are other surfaces which have similar properties and
for which Proposition 2.1 and Lemma 2.2 hold as well.
As is proved in \cite{KK}, $\tilde X$ can be replaced,
in particular, by a fake projective plane.
Recall that, by definition, a {\it fake projective plane} is
a surface of general type with $p_g=q=0$ and $K^2=9.$
Such surfaces do really exist, see \cite{Mu}.

\head 3.\enspace Higher-dimensional examples
\endhead

In this section we prove that
part of Theorem 1.1 which concerns
the dimensions $\ge 3$, the case of surfaces is covered by Section 2.
For this purpose, we develop
several series of examples: the examples
from Proposition
3.7 cover even dimensions $\ge 4$
and those from Proposition 3.8
odd dimensions $\ge 3$.
When applying Propositions
3.7 and 3.8 to prove this part of Theorem 1.1,
it is sufficient to take as $M_1$
the surface constructed
in Section 2 and as $M_2$ a fake projective plane.

Another
series of examples sufficient for the proof
of Theorem
1.1
is given by
Propositions 3.9. and 3.11, see Corollary 3.13.

We start from some auxiliary results which, with one
exception, we did not find
in the literature; they may be of an independent interest.
We use them in the proof of Propositions
3.7, 3.8, 3.9
and 3.11.

\proclaim{Lemma 3.1 (\cite{LS})}
Let $X$ and $X_0$ be compact complex manifolds contained in
a deformation family over
irreducible base space.
Suppose that the canonical class of $X_0$ is ample.
Then, $X$ is a Moishezon variety,
i.e., the transcendence degree of meromorphic function field on $X$
coincides with $\dim X$.
\qed
\endproclaim

\proclaim{Lemma 3.2}
Let $V_1$ and $V_2$ be minimal nonsingular
surfaces of general type.
If there exists a regular map $f:V_1\to V_2$ of degree
$d\geq 1$, then
$K^2_{V_1}\geq K^2_{V_2}$.
In addition,
$K^2_{V_1}= K^2_{V_2}$ if and only if $f$ is an isomorphism
{\rom(}i.e., $d=1$\rom).
\endproclaim

\demo{Proof}
According
to the pull-back formula, $K_{V_1}=f^*(K_{V_2})+R$, where
$R=\sum r_iR_i$ is the ramification divisor of $f$,
$R_i$ are irreducible components of $R$, $r_i\geq 0$
(by ramification divisor we mean here the locus of points
where the jacobian matrix is not of maximal rank).
On the other hand, $K_{V_1}^2>0$, $K_{V_2}^2>0$, and $(f^*(K_{V_2}),R)\geq 0$,
since, by Bombieri theorem (see, f.e., \cite{BPV}),
the linear system $\mid 5K_{V_2}\mid $ has no fixed components.
Besides, $(K_{V_1},R_i)=(f^*(K_{V_2})+R,R_i)\geq 0$
for each irreducible curve
$R_i$ lying on a minimal model.
Therefore,
$$\alignat 2
K_{V_1}^2= & dK_{V_2}^2+2(f^*(K_{V_2}),R)+R^2 \\
= &  dK_{V_2}^2+(f^*(K_{V_2}),R)+(f^*(K_{V_2})+R,R)   \\
= & dK_{V_2}^2+(f^*(K_{V_2}),R)+\sum r_i(f^*(K_{V_2})+R,R_i)
\geq K_{V_2}^2
\endalignat
$$
if $d\geq 1$ and we have the equality iff $d=1$ and $R=0$.
Indeed,
any degree one regular map $f$ is an isomorphism, since
$V_1$ and $V_2$ are minimal surfaces.
\qed
\enddemo

\proclaim{Lemma 3.3}
Let $U$ be a Moishezon variety and $V$ a projective variety having
no rational curves. Then, any meromorphic map $f:U\to V$ is holomorphic.
\endproclaim

\demo{Proof}
By Hironaka theorem,
there is
a sequence $\sigma _i:U_i\to U_{i-1}$, $U_0=U$, $1\leq i\leq k$,
of monoidal transformations with non-singular centers such that
$f_k=f\circ \sigma :U_k\to V$ is holomorphic, where
$\sigma =\sigma _1\circ \dots \circ \sigma _k$.
Note that for any $p\in U$
the preimage $\sigma ^{-1}(p)$ is rationally connected,
i.e., any two distinct points from $\sigma ^{-1}(p)$ belong to a
connected union of a
finite number of rational curves.
Therefore, $f_k(\sigma ^{-1}(p))$ is a point, since
there are no rational curves lying on $V$. Thus, $f_k$ factors through
the holomorphic map $f$. \qed
\enddemo

\proclaim{Lemma 3.4}
Let $Z$ be a compact Moishezon manifold and $X$ be a projective
one. Then, every holomorphic map $f\: Z\to X$ inducing an isomorphism
in homology, $f_*\: H_*(Z;\Z)\to H_*(X;\Z)$, should be biholomorphic.
\endproclaim

\demo{Proof}
To show that
$f$ is an isomorphism, it is sufficient to
check that there is no complex subvariety $Y$ of $Z$ such that
$\dim f(Y)<\dim Y$.

First, let us eliminate the case $\dim Y=\dim Z-1$.
If $\dim f(Y)<\dim Y$ then $Y$ is homological to $0$,
since $f_*$ is an isomorphism. But, since
$X$
is a projective variety, one can find a curve
$B\subset X$ meeting $f(Y)$ and not lying in $f(Y)$.
Therefore, the strict preimage $f^{-1}(B)$ and $Y$ (we assume here that
$Y=f^{-1}(f(Y))$) have positive intersection number, which
contradicts to $Y\simeq 0$.

Now, consider any $Y$ such that $\dim f(Y)<\dim Y$. By above,
we have $\dim Y<\dim Z-1=\dim X-1$. One can find a meromorphic volume
form $\omega$ on $X$
such that $f(Y)$ does not belong to the support of the divisor
$(\omega )$. But it is impossible, since in this case in a neighborhood
of generic point of $Y$ the form
$f^*(\omega)$ would have zero along a submanifold of
$X$ of
codimension at least two. \qed
\enddemo

\proclaim{Lemma 3.5} Let $M_1$ and $M_2$ be nonsingular
regular
surfaces of
general type. If
$M_1$, $\bar M_1$, $M_2$, $\bar M_2$ are pairwise non isomorphic
and contain no rational curves, then:
\roster
\item "{(i)}"
the products $M_1^a\times \bar M_1^b\times M_2^c\times\bar M_2^d$
with $(a,b,c,d)\in\N^4$
are pairwise non isomorphic; in particular, such a
product admits an anti-automorphism only if $a=b$
and $c=d$;
\item "{(ii)}"
in a deformation $\pi\: X\to B_1$
each fiber manifold $X_z=\pi^{-1}(z), z\in B_1$, is isomorphic
to a product $M_1^{j}\times \bar M_1^{m}\times M_2^{k}\times
\bar M_2^{n}$
with $j=j(z), m=m(z),k=k(z), n=n(z)$, only if the deformation is trivial.
\endroster
\endproclaim

\demo{Proof} Let us first prove $(i)$.
Assume that there is an isomorphism
$$h:M_1^{j_1}\times \bar M_1^{m_1}\times M_2^{k_1}\times
\bar M_2^{n_1}\to M_1^{j_2}\times \bar M_1^{m_2}
\times M_2^{k_2}\times \bar M_2^{n_2}.$$
Without loss of generality, we can assume also
that $K_{M_1}^2\geq K_{M_2}^2$ and $j_1\le j_2$.

Denote by $S_l \subset
M_1^{j_1}\times \bar M_1^{m_1}\times M_2^{k_1}\times
\bar M_2^{n_1}$, $l\ge 1$,
a fiber of the projection of the whole product to the
partial product taken over all the factors
except the one with index $l$. So, $S_l$
is canonically isomorphic to $M_1$ for $1\leq l \leq j_1$,
to $\bar M_1$ for $j_1<l\le j_1+m_1$, to $M_2$ for
$j_1+m_1<l\le j_1+m_1+k_1$, and to $\bar M_2$ for $j_1+m_1+k_1<l$.
Denote also by $\Sigma_s$ the similar fibers of
$M_1^{j_2}\times \bar M_1^{m_2}\times M_2^{k_2}\times
\bar M_2^{n_2}.$

Consider, first,
the restrictions $p_{s\mid _l}$ to $h_1(S_l)$ of the projections
$p_s:M_1^{j_2}\times \bar M_1^{m_2}\times M_2^{k_2}\times
\bar M_2^{n_2} \to \Sigma_s$.
Let us show that for each $l\geq j_1+1$
the image
$p_{s\mid _l}(h_1(S_l))$ is a point as soon as $s\leq j_2$.
In fact, $p_{s\mid _l}(h_1(S_l))$ can not be a curve, since
there are no rational curves in $M_1, \bar M_1, M_2$ or in $\bar M_2$, and
if $p_{s\mid _l}(h_1(S_l))$
was a curve of positive genus $g$, then
the irregularity of $S_l$ would be at least $g$.
And since $M_1$
is not isomorphic to any of $\bar M_1, M_2,$ and $ \bar M_2$,
Lemma 3.1 can be applied and it implies that
$p_{s\mid _l}(h_1(S_l))$ can not be a surface neither.

Therefore,
$$h_*:H_2(M_1^{j_1}\times \bar M_1^{m_1}\times M_2^{k_1}\times
\bar
M_2^{n_1})\to H_2(M_1^{j_2}\times \bar M_1^{m_2}\times M_2^{k_2}\times
\bar M_2^{n_2})$$
maps the direct K\"unneth summand
$H_2(\bar M_1^{m_1})\oplus H_2(M_2^{k_1})\oplus
H_2(\bar M_2^{n_1})$
into the direct K\"unneth  summand
$H_2(\bar M_1^{m_2})\oplus H_2(M_2^{k_2})\oplus
H_2(\bar M_2^{n_2})$.
So, $h_*$ can be an isomorphism only if $j_1=j_2$.
Similarly, $m_1=m_2$ (one can simply exchange $M_1$
and $\bar M_1$ in the above arguments).

Now, consider a multi-fiber isomorphic to $M_2^{k_1}\times
\bar M_2^{n_1}.$ From what is proved above it follows that
such a multi-fiber is
mapped by $h$ to a multi-fiber of
$M_1^{j_2}\times \bar M_1^{m_2}\times M_2^{k_2}\times
\bar M_2^{n_2}$
isomorphic to
$M_2^{k_2}\times
\bar M_2^{n_2}$. Applying Lemma 3.2 and the arguments as above
to such a map $M_2^{k_1}\times
\bar M_2^{n_1}\to M_2^{k_2}\times
\bar M_2^{n_2}$, we get, first, that it is an isomorphism and then that
$k_1=k_2$ and $n_1=n_2$.

To prove $(ii)$, let us consider
a deformation $\pi\: X\to B_1$ such that
each fiber manifold $X_z=\pi^{-1}(z), z\in B_1$, is isomorphic
to a product $M_1^{j}\times \bar M_1^{m}\times M_2^{k}\times
\bar M_2^{n}$
with $j=j(z), m=m(z),k=k(z), n=n(z).$
Since neither of the factors contains a rational curve,
all the manifolds $X_z, z\in B_1,$
have ample canonical bundles. Hence, by the Kodaira vanishing
theorem, $\dim H^i(X_z, mK_z)=0$
for any
$m>1$,
$i\ge 1$, and the Riemann-Roch theorem
imples that $\dim H^{0}(X_z, mK_z)$ is constant with respect to $z$.
By Grauert's continuity theorem (see \cite{G}, Section 7, Theorem 4),
it implies that the canonical map $B_1\to\M$,
where $\M$ is the moduli space of projective complex
structures on the smooth oriented manifold underlying $X_z, z\in B_1,$
is defined and holomorphic. Since, under our hypothesis,
the image of this map is countable, it is constant.
\qed
\enddemo

Note,
that when it is only question of applying Lemma 3.5 $(i)$
to products of Miayoka-Yau surfaces,
this Lemma can be replaced by more traditional
arguments based on
the De Rham theorem of uniqueness of the decomposition in irreducible
factors (see, f.e., \cite{KN})
applied to the universal covering equipped with
a Bergmann metric.

\proclaim{Lemma 3.6} Let $M_1$ and $M_2$ be nonsingular
regular surfaces of
general type and $C_1$ and $C_2$ be
non-singular curves of genus $g>0$. If
$M_1, \bar M_1, M_2, \bar M_2$ are pairwise non isomorphic
and contain no rational curves, then the
products $M_{a_1,b_1,c_1,d_1,C_1}=
M_1^{a_1}\times \bar M_1^{b_1}\times M_2^{c_1}\times\bar M_2^{d_1}
\times C_1$ and
$M_{a_2,b_2,c_2,d_2,C_2}=
M_1^{a_2}\times \bar M_1^{b_2}\times M_2^{c_2}\times\bar M_2^{d_2}
\times C_2$, $(a_i,b_i,c_i,d_i)\in\N^4$,
are isomorphic if and only if $(a_1,b_1,c_1,d_1)=(a_2,b_2,c_2,d_2)$
and $C_1$ is isomorphic to $C_2$.
In particular, such a product admits an anti-automorphism only if
$a_1=b_1$ and $c_1=d_1$.
\endproclaim

\demo{Proof}
Since $M_1$ and $M_2$ are regular surfaces,
the image $\alpha _i(M_{a_i,b_i,c_i,d_i,C_i})$
of the Albanese map $\alpha _i:M_{a_i,b_i,c_i,d_i,C_i}\to
\text{Alb}(M_{a_i,b_i,c_i,d_i,C_i})$,
$i=1,2$, coincides with $C_i$.
To conclude it suffices to apply the universal
property of the Albanese map and Lemma 3.2.
\enddemo

\proclaim{Proposition 3.7}
Let
$M_1, \bar M_1, M_2$  and $\bar M_2$
be pairwise non-isomorphic
compact regular nonsingular
surfaces of general type
satisfying the Mostow-Siu
rigidity and having no rational curves.
Let $m$ and $n$ be two non-negative integers and
$X_{0,0}$ be a product
$$S_1\times \dots \times S_m \times S_{m+1}\times \dots \times S_{m+n} ,$$
where each $S_i$, $1\leq i\leq m$, is isomorphic to $M_1$ and
each $S_i$, $m+1\leq i \leq m+n$, is isomorphic to $M_2$.
Then:
\roster
\item "{(i)}" there are at least
$(m+1)(n+1)$ distinct deformation classes of complex structures on
the underlying oriented smooth manifold $X_{0,0}$; these classes
are represented by
$$X_{m-j,n-k}=
M_1^j\times\bar M_1^{m-j}\times
M_2^{k}\times \bar M_2^{n-k},$$
with $ j=0,\dots , m$ and $k=0,\dots , n$
(here, $M_l^j$ states for the product of $j$ copies of $M_l$
and $\bar M_l^s$ for the product of $s$ copies of $\bar M_l$);
\item "{(ii)}" there are no anti-automorphisms on $X_{p,q}$
except the case $m=2p$, $n=2q$;
for other values of $p,q$ the deformation class
of $X_{p,q}$ is not invariant
under reversing of complex structure, $\bar X_{p,q}=X_{m-p,m-q}$.
\item "{(iii)}"
a complex manifold is deformation equivalent to $X_{p,q}$
if and only if it is isomorphic to $X_{p,q}$.
\endroster
\endproclaim

Note that, as it follows from the proof
below, the complex structures represented by $X_{p,q}$
are the only Moishezon complex structures on the smooth oriented
manifold underlying $X_{0,0}$.
Up to our knowledge, the existence of non
Moishezon complex structures on this manifold
is an open question. If all the complex structures on $X_{0,0}$
are Moishezon, the number of their deformation classes is exactly
$(m+1)(n+1)$.

\demo{Proof}
Clearly, $(i)$ and $(ii)$
are straightforward consequences of $(iii)$. So, let us prove the part $(iii)$
of the statement.

As it follows from Lemma 3.1 and  Lemma 3.5 $(ii)$,
it is sufficient to show that
any Moishezon complex structure on the underlying $X_{0,0}$ smooth
manifold provides a complex manifold isomorphic to one of $X_{p,q}$.

So, consider a Moishezon complex compact variety
$X$ having the same underlying smooth
oriented manifold as $X_{0,0}$.
Pick a diffeomorphism between $X$ and $X_{0,0}$ and consider
the corresponding to it projections $p_i:X\to S_i$,
$i=1,\dots ,m+n$, (we use here notations
from Lemma 3.5).
By Moishezon Theorem, there is
a sequence $\sigma _i:Z_i\to Z_{i-1}$, $Z_0=X$, $1\leq i\leq r$,
of monoidal transformations with non-singular centers such that
$Z=Z_r$ is a projective variety. Denote by $\sigma$ their composition.
Every
projection and every monoidal transformation
induce an epimorphism in homology.
So, according to the
Mostow-Siu rigidity hypothesis,
each $p_i\circ\sigma$, $i=1,\dots ,m+n$, is
homotopic to a map
$\widetilde{p}_i:Z\to S_i$
which is either holomorphic or anti-holomorphic. Let $j$
be the number of holomorphic maps $\widetilde{p}_i$ for $i\leq m$ and
$k$ be the number of holomorphic maps $\widetilde{p}_i$ for $i> m$.
Thus, after a suitable renumbering of $\widetilde p_i$,
we get a meromorphic map
$f=\widetilde{p}_1\circ\sigma^{-1}\times \dots \times \widetilde{p}_{m+n}
\circ\sigma^{-1}
\: X \to M_1^j\times \bar M_1^{m-j}\times
M_2^{k}\times \bar M_2^{n-k}$.
By Lemma 3.3 it is holomorphic, and, since $\sigma$
induces an epimorphism in homology,  from the
construction of $f$ it follows that $f$ induces an isomorphism in homology.
Now, Lemma 3.4 applies and it shows that $X$ is isomorphic to $X_{j,k}$.
\qed
\enddemo

\proclaim{Proposition 3.8}
Let $M_1$ and $M_2$ be as in Proposition 3.7,
$C_0$ be a curve of genus $g>1$, and
$Y_{0,0}$ be a product $M_1^m \times M_2^{n}\times C_0$.
Then:
\roster
\item "{(i)}" there are at least $(m+1)(n+1)$
distinct deformation classes of complex structures
on the underlying oriented smooth manifold $Y_{0,0}$;
these deformation classes are represented by
$$Y_{m-j,n-k}=
M_1^j\times \bar M_1^{m-j}\times
M_2^{k}\times \bar M_2^{n-k}\times C_0,$$
$j=0,\dots , m$ and $k=0,\dots , n$;
\item "{(ii)}" there
are no anti-automorphisms on $Y_{p,q}$
except the cases when $m=2p$, $n=2q$
and $C_0$ admits an anti-automorphism;
for all the other values of $p,q$
the deformation class of $Y_{p,q}$ is not invariant
under reversing of complex structure, $\bar Y_{p,q}=Y_{m-p,m-q}$.
\endroster
\endproclaim

Here, the dimension of complex manifolds $Y_{p,q}$ is odd,
and, thus, reversing the complex structure
in 3.8(ii) we change the complex orientation.
However, since $Y_{p,q}$ admit orientation reversing diffeomorphisms,
there are as much deformation classes of complex structures on
the underlying smooth manifold
as with fixed or with any orientation.

Similarly to the situation of Proposition
3.7, if all the complex structures on $Y_{0,0}$
are Moishezon, the number of their deformation
classes is exactly $(m+1)(n+1)$.

\demo{Proof}
Let start of the first part of Proposition.
As in the proof of Proposition 3.7,
to prove the first part
we should check that the complex structure
which can be deformed to a Moishezon complex structure
is isomorphic to some
$M_1^j\times \bar M_1^{m-j}\times
M_2^{k}\times \bar M_2^{n-k}\times C$, where $C$ is a curve
of genus $g$, and then we should prove that  a complex structure
which is deformation equivalent to one of $Y_{p,q}$
is isomorphic to the complex structure on $M_1^p\times \bar M_1^{m-p}\times
M_2^{q}\times \bar M_2^{n-q}\times C$ for some curve $C$.

First, consider a Moishezon variety
$Y$ having the same underlying smooth
oriented manifold as $Y_{0,0}$.
Consider the Albanese map $\alpha : Y \to \text{Alb}(Y)$.
The Albanese map
$\alpha_{0,0} : Y_{0,0} \to \alpha_{0,0}(Y_{0,0})\subset
\text{Alb}(Y_{0,0})$ coincides with
the last factor projection $p_{m+n+1} :Y_{0,0}=M_1^m \times
M_2^{n}\times C_0\to C_0$.
Since the irregularity
$q(M_1^m \times M_2^{n})$
is zero, the subring of $H^*(Y_{0,0};\Z)$
generated by $H^1(Y_{0,0}; \Z)$
coincides with $\alpha_{0,0}^{*}(H^{*}(C_0;\Z))$.
In particular, $\bigwedge^{2} H^1(Y_{0,0}, \Z)\subset H^{2}(Y_{0,0},\Z)$
is a one-dimensional subspace and
$\bigwedge^{i} H^1(Y_{0,0}, \Z)=
\{ 0\}$ for $i>2$. Therefore, the image $\alpha (Y)$ is a curve $C$ of
genus $g$ and $\alpha^{*}(H^{*}(\text{Alb}(Y),\Z))$ also
coincides with the subring generated by $H^1(Y,\Z)=H^1(Y_{0,0}, \Z)$.

Pick a diffeomorphism between $Y$ and $Y_{0,0}$ and consider
the corresponding to it projections $p_i:Y\to S_i$,
$i=1,\dots ,m+n$, (we use here notations
from Proposition 3.5), which are $C^\infty$-submersions,
and $p_{m+n+1}:Y\to C$.
By Moisheson Theorem, there is
a sequence $\sigma _i:Z_i\to Z_{i-1}$, $Z_0=Y$, $1\leq i\leq k$,
of monoidal transformations with non-singular centers such that
$Z=Z_k$ is a projective variety. Denote by $\sigma$ their composition.
Every
projection and every monoidal transformation
induce an epimorphism in homology.
So, according to the
Mostow-Siu rigidity hypothesis,
each $p_i\circ\sigma$, $i=1,\dots ,m+n$, is
homotopic to a map $\widetilde{p}_i:Z\to S_i$
which is either holomorphic or anti-holomorphic. Let $j$
be the number of holomorphic maps $\widetilde{p}_i$ for $i\leq m$ and
$k$ be the number of holomorphic maps $\widetilde{p}_i$ for $i> m$.
Thus, we get a meromorphic map
$f=\widetilde{p}_1\circ\sigma^{-1}\times \dots \times \widetilde{p}_{m+n}
\circ\sigma^{-1}
\times \alpha : Y \to M_1^j\times \bar M_1^{m-j}\times
M_2^{k}\times \bar M_2^{n-k}\times \alpha(Y)$.
By Lemma 3.3 it is holomorphic, and, since $\sigma$
induces an epimorphism in homology,  from the
construction of $f$ it follows that $f$ induces an isomorphism in homology.
Now, Lemma 3.4 applies and it shows that $Y$ is isomorphic to $Y_{j,k}$.

To finish the proof of the first part of
Proposition, let us consider
a deformation $X\to B_1$ of $Y_{m-j,n-k}.$
As it is already proved, for any $z\in B_1$
the manifold $X_z$ is isomorphic to
$M_1^j\times \bar M_1^{m-j}\times
M_2^{k}\times \bar M_2^{n-k}\times C_z$
with some $j=j(z), k=k(z)$.
Hence, as it follows, for example, from
Grauert's continuity theorem (see \cite{G})
the associated Albanese
varieties $\Alb X_z$ of $X_z, z\in B_1,$
form a deformation $W\to B_1$.
The Albanese map $X\to W$ is of constant rank equal to
$2$ (one unit comes from the images $C_z\subset\Alb X_z$ of the
Albanese map and another one from $B_1$, since
as any deformation the composite
map $X\to W\to B_1$ is a submersion). It is proper
and defines a deformation of
$M_1^{j}\times \bar M_1^{m-j}\times
M_2^{k}\times \bar M_2^{n-k}$.
By Proposition 3.7 (iii),
all elements of such a deformation are isomorphic to each other.

The second part of Proposition is a straightforward consequence of the
first part and Lemma 3.6.
\qed
\enddemo

{\bf Remark.} If one restricts
himself to K\"ahler (or Moishezon) manifolds
and there deformations constituted only of K\"ahler
(or Moishezon) manifolds,
then using the same arguments as in the proof of Propositions
3.7 and 3.8 one can prove
that $X_{p,q}\times T_1$ and $X_{k,l}\times T_2$,
where $T_1, T_2$ are tori of equal dimensions and
$X_{p,q}, X_{k,l}$ are as in Proposition 3.7,
are equivalent if and only if $p=q$ and $k=l$.
In particular,
$\bar X_{p,q}\times T_1$ and $X_{p,q}\times T_2$
are not equivalent except the case $2p=m, 2q=n$
(for the definition of $m$ and $n$ see Proposition 3.7).

\proclaim{Proposition 3.9} Let $M$ be a fake projective plane or
the surface constructed in section 2, and let
$X$ be a K\"ahler
manifold of Kodaira dimension $\leq 1$ with $c_1(X)\in H^2(X,\Z )$
divisible by some
integer $p>3, p\neq 6$ (in particular,
$c_1(X)$ can be zero). Then, there is no
homeomorphism
$h : Y\to Y$,
$Y=M\times X,$ such that $h^*(c_1(Y))=-c_1(Y)$.
\endproclaim
\demo{Proof}
Assume that there is such a
homeomorphism $h$.
Since $M$ is a regular surface,
the K\"unneth decomposition of $H^2(Y;\Z)$
takes the form $H^2(Y,\Z)=H^2(M,\Z )\oplus H^2(X,\Z )$.
Let $h^*(c_1(M))=\omega_1 +\nu _1$ and
$h^*(c_1(X))=\omega_2+\nu _2$, where $\omega _i\in H^2(M,\Z )$ and
$\nu _i\in H^2(X,\Z)$. Then, according to our assumption,
$\omega_1+\omega_2=-c_1(M)$ and $\nu_1+\nu_2=-c_1(X)$.

Denote by
$p_1:Y\to M$
the projection to the first factor,
by $p_2:Y\to X$ the
projection to the second one,
by $i_2: X\to Y$
the canonical isomorphism with
one of the fibers of $p_1$,
and by $i_1:M\to Y$
the canonical isomorphism with
one of the fibers of $p_2$.
Consider the
composition $p_1\circ 
h$.
By Mostow-Siu
rigidity, $p_1\circ 
h$ is homotopic
to $\widetilde p_1:
Y
\to M$ which is either holomorphic or anti-holomorphic
dominant map.

Since $M$ is regular and has no rational curves,
its image $(\tilde p_1\circ i_1) M\subset M$
can not be a curve,
and
it can not be a point, since in this case the
restriction of $\widetilde p_1$ to $i_2(X)$
should be a dominant map, but
the latter is
impossible because of
$\text{Kod dim}X< \text{Kod dim}M$. Therefore,
Lemma 3.2 applies
and it shows that $\widetilde p_1\circ i_1$ is a
holomorphic
or anti-holomorphic isomorphism. Since $M$ has no anti-automorphisms,
$\widetilde p_1\circ i_1$
is a holomorphic map.
Thus, $\omega_1=(p_1\circ h\circ i_1)^* c_1(M)=(\tilde p_1\circ i_1)^*
c_1(M)=c_1(M)$.
It implies
$(p_1\circ h\circ i_2)^*c_1(X)=\omega _2= -2c_1(M)$,
and to complete the proof
it remains to
note that $2c_1(M)$ is not divisible by
$p>3, p\neq 6$,
which follows from
$c_1^2(M)=9$ in the case of
fake projective planes
and from $c_1^2(M)=333$
in the case of the surface constructed in Section 2.\qed
\enddemo

\proclaim{Corollary 3.10}
Let
$Y$
be a variety
as in Proposition 3.9.
Then,
it is not deformation
equivalent to its conjugate and, in particular,
can not be deformed to a variety with anti-automorphisms.\qed
\endproclaim

{\bf Remark.}
In the proof of Proposition 3.9 the hypothesis
on Kodaira dimension of $X$ is used only to exclude existence
of dominant maps from $X$ to $M$.

One more series of manifolds whose canonical class and its inverse
are topologically distinct is given by the following proposition.

\proclaim{Proposition 3.11} Let $M_1$ be the surface constructed in section 2,
$M_2$ a fake projective plane, $X$ as in Proposition 3.9, and
$Y=M_1\times M_2^a\times X$. Then there is no homeomorphism $h: Y\to Y$
such that $h^*c_1(Y)=-c_1(Y)$.
\endproclaim

The proof of Proposition 3.11 is
essentially the same as the proof of Proposition 3.9.
The principal new element with respect to that proof
(except the analysis of the first component of the
homeomorphism $h$, where we use arguments
like in the proof of Proposition 3.7 in order to prove that the
restriction to $M_1$ of the corresponding $\widetilde p_1$ is an isomorphism)
is the following lemma applied to $Z=M_2$.

\proclaim{Lemma 3.12} Let
$Z$ be a surface of general type such that for any
holomorphic map
$f:Z\to Z$ either $\codim f(Z)\geq 2$ or $f$ is a
biholomorphic. Then, the same is true for holomorphic maps
from $Z^a$ to $Z^a$ for any $a\ge 1$.
\endproclaim

\demo{Proof}
Let $f:Z^a\to Z^a$ be a holomorphic map. Denote by
$p:Z^a\to Z^{a-1}$
one of the canonical projections and by
$i:Z\to Z^a$
one of its fibers. If the projections of
$(f\circ i)Z$
to each of the factors are all not dominant, the image
of this fiber is a point. Then, the image of all the parallel
fibers is a point, which implies $\codim f(Z^a)\geq 2$.

Otherwise, one of these projections provides an isomorphism
$Z\to Z$.
Since the automorphism group of $Z$ is discrete,
the map $f$ has, up to permutation of the factors in the source and
the target, the following triangle form
$f(x,p)=(gx, h(x,p))$, where
$g$ is an automorphism of $Z$ and $h$ is a holomorphic
map $Z^a\to Z^{a-1}$. By induction, we may assume that
for some, and hence for any, $x\in Z$ the map $Z^{a-1}\to Z^{a-1}$
given by
$p\mapsto h(x,p)$ either is biholomorphic or has a codimension
$\geq 2$ image. The above triangle form
of $f$ implies that it has the same property.
\qed
\enddemo

\proclaim{Corollary 3.13} For any $n\geq 2$ there is a projective
$n$-dimensional manifold $Y$ of general type which has no homeomorphisms $h$
with $h^*c_1(Y)=-c_1(Y)$.
\endproclaim

\demo{Proof} Apply Proposition 3.11 to a curve $X$ of genus $g>4$
in the odd-dimensional case and to a point $X$ in the even-dimensional case.
\qed
\enddemo

In fact,
the same arguments which were used
in the proof of Proposition 3.7 (or more traditional arguments
base on the De Rham theorem of uniqueness of the
decomposition in irreducible factors)
allow
to deduce from Lemma 2.3 a complete description
of the homeotopy group $\Cal{H}(Y)$
in the case $Y=M^n$, where $M$ is the surface
constructed in Section 2.

\proclaim{Proposition 3.14} Let $M$ be the surface constructed in
Section 2. Then, the  group $\Cal{H}(M^n)=
(\Z /5\Z \oplus \Z /5\Z)^n\rtimes S_n$ is the semi-direct product, with
standart action
of the symmetric group $S_n$ on the factors $(\Z /5\Z \oplus \Z/5\Z )$ of
the direct product $(\Z /5\Z \oplus \Z /5\Z)^n$.\qed
\endproclaim

\head 4.\enspace Applications
\endhead

{\bf A. Non connected moduli spaces.}
Here,
we discuss how
Theorem 1.1 and
our solution of the $\dif=\Def$ problem
can be translated
in the language of moduli spaces.

In
most of our examples
the 
varieties 
are of general type
and, moreover, all
varieties deformation equivalent
to them are
projective varieties
of general type
isomorphic to products of projective surfaces of general type
and, possiblly, some curve of
genus $>1$
(see Propositions 2.1, 3.7, and 3.8).
So,
here
we may
restrict ourselves to the case when
the {\it moduli space $\M_{\dif}$ of complex structures} on
a given smooth orientable
manifold $M$ is well defined as a complex
space.

In fact, we fix the orientation of $M$
only if the complex dimension of $M$
is even. Then, obviously, in even
complex dimensions
the $\dif=\Def$ problem becomes the question
of connectedness of $\M_{\dif}$.
In addition, if the complex dimension is even,
the involution $\c: X=(M,J)\mapsto\bar X=(M,-J)$
acts on $\M_{\dif}$ and defines a canonical real structure on it.

Since
in all our odd-dimensional examples
$M$ has orientation reversing diffeomorphisms,
the $\dif=\Def$ problem as it is stated in Introduction
remains
the question
of connectedness of $\M_{\dif}$
even without introducing orientations
in the definition of the moduli space. On the other hand,
it allows still to consider the complex conjugation
acting on $M_{\dif}$ by $\c: X=(M,J)\mapsto\bar X=(M,-J)$.

To reinterpret Theorem 1.1
in its full content is more complicated,
in general. Certainly,
if $X$ can be equipped with a real structure then it defines a real point
of $\M_{\dif}$, i.e., a fixed point of $\c$. But the reverse is not true, in
general, since real points of $\M_{\dif}$
are given also by complex manifolds
having antiholomorphic automorphisms of any (even) order.
The simplest examples are given by the Shimura curves \cite{Shi}
which are double coverings of $\Cp1$ ramified in $4m+2\ge 6$ generic points
invariant under the real structure of $\Cp1$ acting without fixed points.
The real structure is lifted
to an anti-automorphism $\c$ of order $4$, and for
generic points the curve has no other anti-automorphisms than $\c$ and
the other lift, which is its composition with a deck transformation and which
is also of order $4$. In fact, in these examples the automorphism group
$\Aut$ is $\Z/2$ and the Klein group $\Kl$, which includes all the
anti-automorphisms and all the automorphisms, is $\Z/4$.

By contrary, the absence of deformation equivalence between $X$ and $\bar X$
has a simple meaning in terms of $\M_{\dif}$: it means that
$\M_{\dif}$ contains two connected components interchanged by $\c$.
Thus, the following statement is a straightforward
corollary of Theorem 1.1 and Propositions 3.7, 3.8.
(Let us recall our convention: in this
theorem, as well as above, we fix the orientation only if
the complex dimension of the manifold
is even.)

\proclaim{Theorem 4.1}%.\enspace Theorem}
In any even real dimension $2n\ge 4$ there exist an oriented smooth
compact manifold for which the moduli space of complex structures
is disconnected and has at least
$([\frac{n}4]+1)([\frac{n}2]-[\frac{n}4]+1)=N$
connected components. At most
one of the components is invariant under complex conjugation.
\endproclaim

This statement can be refined a bit.
For example, in the
case of surfaces, the moduli space
consists of two connected components
and the complex structures
from different components have opposite canonical classes.

{\bf B. Isotopy classes of cuspidal curves.}

The theorem below shows that diffeomorphic plane cuspidal curves can be
not equivalent under equisingular deformations. In fact, in our examples the
curves are not even isotopic. Thus, the following question remains open:
does existence of an isotopy imply existence of an equisingular deformation.

\proclaim{Theorem 4.3} There are two infinite
sequences, $\{C_{m,1}\}$ and $\{C_{m,2}\}$,
of plane irreducible cuspidal curves
of degree $\deg(C_{m,1})=\deg(C_{m,2})\to\infty$,
such that the pairs $(\Cp2,C_{m,1})$ and $(\Cp2,C_{m,2})$
are diffeomorphic, but $C_{m,1}$ and $C_{m,2}$ are not isotopic and,
in particular, they can not be connected by an equisingular deformation.
\endproclaim

To get Theorem 4.3 it is sufficient to apply Proposition 4.4 stated below
to any of the surfaces given in Section 2.

\proclaim{Proposition 4.4}
Let $X$ be a surface of general type with ample
canonical class $K$. Suppose that there is no
homeomorphism $h$ of $X$ such that $h^*[K]=-[K]$,
$[K]\in H^2(X;\Q)$.
Then the moduli space of $X$ consists of at least 2
connected components corresponding to
$X$ and $\overline{X}$ (the bar states for reversing of complex structure,
$J \mapsto -J$), and for any $m\geq 5$
these two connected components are distinguished by
the isotopy types of the branch curves of generic coverings
$f_m:X\to \Cp2$ and $\overline{f}_m:\overline{X}\to \Cp2$
given by $mK$ and $m\overline K$, respectively. In particular,
the branch curves can not be connected by an equisingular deformation.
\endproclaim

When Proposition 4.4 is applied to the surface $X=\tilde X$
constructed in Section 2 by means of the C\'eva configuration,
one finds (see general formulae in
\cite{K}
, page 1155)
\roster
\item"{(1)}" $\deg f_m=333\,m^2$;
\item"{(2)}" $C_{m,1}$ is a plane cuspidal curve of
$\deg C_{m,1}= 333\, m(3m+1)$;
\item"{(3)}" the geometric genus $g_m$ of $C_{m,1}$ equals
$g_m=333\,(3m+2)(3m+1)/2+1$;
\item"{(4)}" $C_{m,1}$ has $c=111\,(36m^2+27m+5)$ ordinary cusps.
\endroster
If it is a fake projective plane which is taken as $X$, then
$$
\deg f_m=9\,m^2,\quad
\deg C_{m,1}= 9\, m(3m+1),
$$
$$
g_m=\frac92 \,(3m+2)(3m+1)+1, \quad
c=3\,(36m^2+27m+5).
$$

{\it Sketch of the proof of Proposition 4.4} \rom(see \cite{KK2} for
more details). Since $K$ is ample, from Bombieri theorem
it follows that the map $X\to\C p^{r_m}$,
$r_m=\dim H^0(X,mK)-1$, given by $mK$ is an imbedding
if $m\ge 5$. Let $m\ge 5$ and denote
by $X_m$ the image of $X$ under
the imbedding in $\C p^{r_m}$ given by $mK$, by
$pr_m: \C p^{r_m}\to \Cp2$ a linear projection generic with respect to
$X_m$, by $f_m=pr_m|_{X_m}$ the restriction of $pr_m$ to $X_m$, and
by $C_{m,1}\subset \Cp2$ the branch curve of $f_m$.
As soon as we identify  $X_m$ and $\overline{X}_m$ as sets,
the composition
$\overline{f}_m=c\circ f_m :\overline{X}_m\to \Cp2$
of $f_m$ with the standard complex conjugation
$c:\Cp2\to\Cp2$ is a
holomorphic generic covering with branch curve
$C_{m,2}=c(C_{m,1})$. By construction, we have
$$\overline{f}_m^{\,*}(\Lambda )= -f_m^{*}(\Lambda )=-m[K]\,  ,
\tag 3
$$
where $\Lambda \in H^2(\Cp2,\Q)$ is the class of the projective line in
$\Cp2$.

The set of generic coverings $f$ of $\Cp2$ branched along a cuspidal curve
$C$ is in one-to-one correspondence with the set of epimorphisms from
the fundamental group $\pi _1(\Cp2\setminus C)$ to the symmetric groups
$S_{\deg f}$ (up to inner automorphisms) satisfying some additional
properties (see
%\cite{\refK}
\cite{K}).
By Theorem 3 in \cite{K}, for $C_{m,1}$ (respectively, $C_{m,2}$)
there exists only one such an epimorphism
$\varphi_m: \pi _1(\Cp2\setminus C_{m,1})\to S_{\deg f_m}$
(respectively, $\overline{\varphi}_m$).
Thus, if there exists an isotopy
$F_t: \Cp2 \to \Cp2$ such that $F_0=\id$ and
$F_1 (C_{m,1})=C_{m,2}$, then the epimorphism
$$
\varphi_{m,t}: \pi _1(\Cp2\times [0,1]\setminus
\{ (F_t(C_{m,1}),t)\} )\to S_{\deg f_m}
$$
defines a locally trivial, with respect to $t$, family of generic coverings
$f_{m,t}:Y\to \Cp2, t\in [0,1]$, which provides in its turn
a homeomorphism $h$ of $X$
such that $h^*[K]=-[K]\in H_*(X;\Q)$, as it follows from $(3)$.
Thus, the isotopy $F_t$ can not exist.
\qed

{\bf C. Two inequivalent classes of symplectic structures.}
Here,
we treat the equivalence relation 
defined in Introduction.

Let $Y$ be
a K\"ahler surface
satisfying Lemma 2.2,
or, more generaly, a
K\"ahler manifold $Y$ which has no diffeomorphisms
$f:Y\to Y$ with $f^*[K]=-[K], [K]\in H^2(Y;\Q)$.
One can take, for
example, the surfaces given in Section 2
or
any manifold $Y$ like in
Propositions 3.9
or 3.11.
Denote by $\omega$ the symplectic structure on $Y$ which is
the imaginary part of its K\"ahler structure.

\proclaim{Proposition 4.5}
The symplectic structures
$\omega$ and $-\omega$ are not equivalent to each other.
In particular, they are not symplectomorphic.
\endproclaim

\demo{Proof}
The class $[K]$ is
the canonical class of $\omega$,
while $-[K]$ is the canonical class of $-\omega$.
So, the result follows from
the invariance
of the canonical class under deformations
and the absence of diffeomorphisms
transforming $[K]$ in $-[K]$.\qed
\enddemo

{\bf Remark.} The above argument can be applied to
deformations
in
the class of
almost-complex structures.
It shows that for any
manifold
$Y=(Y,J)$ as above the structures $J$ and $-J$
are not equivalent, where similar to equivalence of symplectic structures,
we call two almost-complex structures {\it equivalent} if they
can be obtained one from
another by deformation followed, if necessary, by a diffeomorphism.

The following proposition shows
that
between manifolds studied in Section 3
there are
more inequivalent
symplectic and almost-complex structures.
Its proof repeats word by word the proof of Proposition 3.7 (iii).
\proclaim{Proposition 4.6}
Let
$Y$ be the common underlying oriented smooth manifold
of complex manifolds $X_{p,q}$ from Proposition 3.7.
Then,
\roster
\item"{(i)}" there are no homeomorphisms $h:Y\to Y$
transforming $c_1(X_{p,q})$ in $c_1(X_{s,r})$ except if
$p=s$ and $q=r$; %. Thus, their
\item"{(ii)}" symplectic
\rom(respectively, almost-complex\rom)
structures on $X_{p,q}$ and $X_{r,s}$ are non equivalent
except if $p=s$ and $q=r$.
\qed
\endroster
\endproclaim

\head Appendix.\enspace Geometric genus calculation
\endhead

The aim of this Appendix is to
show that the irregularity of the surface constructed in Section 2
equals zero.
We start from explaining a general algorithm we use for this
calculation and only after that apply it to the particular case
in question. In fact, we calculate instead the geometric genus, which
is sufficient, since their difference is a topological invariant,
due to Noether's formula. In the calculation we use permanently
the invariance of the geometric genus under birational transformations,
which allows us at each step to use that
nonsingular birational model which is
more convenient for calculation.

The result is by no means new. It is contained, for example, in
the results stated in \cite{I}. For completeness
of the proof of Theorem 4.1,
we present a straightforward
calculation and with as much details as possible.

{\bf Algorithm of reduction to cyclic coverings.}
Let $f:X_G\to \Cp2$, where
$X_G$ is supposed to be a normal
surface, be a Galois covering with
abelian Galois group
$G=(\Z /p\Z)^m$,
where $p$ is a prime number, and
branched along curves $B_1,\dots , B_t\subset \Cp2 $.
Such a covering is determined by
an epimorphism $\phi :H_1(\Cp2 \setminus \cup B_i) \to G$.
Write
it in a form
$$\phi (\gamma _i)=
k_{1,i}\alpha_1+\dots +k_{m,i}\alpha_m,
\qquad i=1,\dots ,\, t,$$
where
$\alpha_j$
are standard generators of $G=\oplus(\Z/p\Z)^m$,
$\gamma _i$ are standard generators of $H_1(\Cp2 \setminus \cup B_i)$
dual to $B_i$ and $k_{i,j}\in \Z /\Z p$,
$0\leq k_{i,j} < p$, are coordinates of
$\phi (\gamma _{i})$ with respect to $\alpha_j$.
In this notation,
$X_G$
is the normalization of the projective closure of
the affine
surface
$Y_G\subset \C^{m+2}$ given by
$$ z_j^{p}=\prod_{i=1}^t h_i^{k_{j,i}}(x,y), \qquad , \, j=1,\dots ,\, m,
$$
where $h_i(x,y)$ are
the equations of $B_i$
in some chart $\C^2 \subset \Cp2$.

Let $\widetilde X_G$ be the minimal desingularization of $X_G$.
As is known,
it exists,
it is unique and
the action of $G$ lifts, in an unique way, to
a regular action on $\widetilde X_G$.

Consider
the action of $G$
on the space $H^0(\widetilde X_G,
\Omega ^2_{\widetilde X_G})$ of regular 2-forms.
It provides a decomposition
$$H^0(\widetilde X_G,
\Omega ^2_{\widetilde X_G})
=\oplus H_{(s_1,,\dots ,s_m)}$$
into the direct sum of eighenspaces $H_{(s_1,\dots, s_m)}$,
where
$$\omega \in H_{(s_1,\dots, s_m)}
\quad\text{iff}\quad
\alpha _j(\omega )=e^{2\pi s_j i/p}\cdot \omega
\quad \text{for any}
\, j=1,\dots , \, m.$$
Let
$H\subset G$ be a subgroup  and $G_1=G/H$.
We have the following commutative diagram
$$
\CD
X_{G}    @>>{h}>                 X_{G_1}   \\
@V{f}VV                    @V{f_1}VV   \\
\Cp2              @>>{id}>           \Cp2
\endCD
$$
where $f_1: X_{G_1}\to \Cp2$
is the Galois covering
corresponding
to $\phi _1=i\circ \phi $
with $i:G\to G_1=G/H$
being the
canonical epimorphism.
The
map $h$ induces a
rational dominant (i.e., whose image is everywhere dense) map
$\widetilde X_{G}\to \widetilde X_{G_1}$,
and the latter, as any
rational dominant map between
nonsingular varieties, transforms holomorphic
$p$-forms in holomoprhic $p$-forms. Thus,
the
subspace $h^*(H^0(\widetilde X_{G_1},\Omega^2_{\widetilde X_{G_1}}))
\subset H^0(\widetilde X_G,\Omega^2_{\widetilde X_G})$
is well defined, and it coincides with
the subspace $H^0(\widetilde X_G,\Omega^2_{\widetilde X_G})^H
\subset H^0(\widetilde X_G,\Omega^2_{\widetilde X_G})$
of the elements
invariant under the action of $H$.
On the other hand, an eighenspace
$H_{(s_1,\dots, s_m)}$ is invariant under
$x_1\alpha_1+\dots x_m\alpha_m$ if and only if $x_1s_1+\dots+x_ms_m=0\, (p)$.
Hence, the sum
$\oplus H_{(\theta s_1,\dots, \theta s_m)}$
taken over $\theta\in\Z/p\Z$
coincides with
$H^0(\widetilde X_G,\Omega^2_{\widetilde X_G})^H$,
where
$$H=\{ \, x_1\alpha_1+\dots +x_m\alpha_m\, \, \vert \, \,
x_1s_1+\dots x_ms_m=0\, (p)\} .$$
So,
this sum is isomorphic to
$H^0(\widetilde X_{G/H},\Omega^2_{\widetilde X_{G/H}}).$
These considerations give rise to the following result.

\proclaim{Proposition 5.1} The geometric
genus $p_g(\widetilde X_G)=\dim H^0(\widetilde X_G,\Omega^2_{\widetilde X_G})$
of $\widetilde X_G$ is equal to
$$ p_g(\widetilde X_G)= \sum_{H} p_g(\widetilde X_{G/H}),$$
where the sum is taken over all subgroups $H$ of $G$ of $rk\, H=rk\, G-1.$
\endproclaim

{\bf Cyclic coverings.}
Now, let $G=\Z /\Z p$ be a cyclic group. To compute
$p_g(\widetilde X_G)$, let us choose homogeneous coordinates
$(x_0:x_1:x_2)$
in $\Cp2$ such that the line $x_0 =0 $
does not belong to the branch locus of $f: X_G\to \Cp2$.
As above, $X_G$ is the normalization of the projective
closure of the hypersurface in $\C^3$ given by equation
$$z^p=h(x,y),$$
where $x=\frac{x_1}{x_0}$,
$y=\frac{x_2}{x_0}$,
$$h(x,y)=\prod_{i=1}^t
h_i^{k_i}(x,y),$$
$h_i(x,y)$ are irreducible equations in $\C^2\subset \Cp2$
of the curves $B_i$
constituting
the branch locus,
and $0<k_i<p$.
Note that the degree
$$
\deg h(x,y)=\sum k_i\deg h_i(x,y)=np
$$
is divisible by $p$,
since the line $x_0 =0 $
does not belong to the branch locus.

It is easy to see that over the chart $x_1\neq 0$ the variety $X_G$
coincides with the normalization of the hypersurface in $\C ^3$ given by
equation
$$w^p=\widetilde h(u,v),$$
where $u=\frac{1}{x}$, $v=\frac{y}{x}$, $\widetilde h(u,v)=u^{np}h(
\frac{1}{u},\frac{v}{u})$, and $w=zu^n$.

{\it Regularity condition over a
generic point of the base.}
Consider $$\omega \in H^0(X_G\setminus \text{Sing} X_G,
\Omega ^2_{X_G\setminus \text{Sing}X_G})$$
and find a criterion of its regularity outside the ramification and
singular loci.

Over the chart $x_0\neq 0 $
the form $\omega$
can be written as
$$\omega =(\sum_{j=0}^{p-1}z^jg_j(x,y))\frac{dx\wedge dy}{z^{p-1}},\tag{4} $$
where $g_j(x,y)$ are rational functions in $x$ and $y$.
The form
$$\frac{dx\wedge dy}{z^{p-1}}$$
has neither poles or
zeros outside of the preimage of the branch locus.
Therefore,
$\omega$ is regular at such a point iff
all $g_j(x,y)$ are regular
at each point
$(a,b)\not\in \sum B_i$.

In fact, if some
$g_j(x,y)$ is not regular at $(a,b)$, then the sum
$$\sum_{j=0}^{p-1}z^jg_j(x,y)$$ can be written
as
$$\frac{\sum_{j=0}^{p-1}z^jP_j(x,y)}{P_p(x,y)},$$
where $P_j(x,y)$, $j=0,\dots ,p$, are polynomials such that
$P_j(a,b)\neq 0$ for some $j<p$
and $P_p(a,b)=0$. Therefore,
$$\sum_{j=0}^{p-1}z^jP_j(a,b)=0$$
at all the $p$ points belonging to
$f^{-1}(a,b)$, since otherwise
$\omega$ would not be regular over $(a,b)$.
On the other hand, it is impossible,
since a non-trivial polynomial of degree less than $p$ can not have
$p$ roots.

{\it Regularity condition over the line at infinity.}
Consider the form $\omega$ over the chart $x_1\neq 0$,
$$\omega = -(\sum_{j=0}^{p-1} w^j\frac{\widetilde g_j(u,v)}{u^{jn
+\deg g_j}})
\frac{1}{u^{3-n(p-1)}}
\frac{du\wedge dv}{w^{p-1}},$$
The similar arguments as above show that
the regularity criterion is equivalent to the following
bound on the degree of the rational functions $g_j$
$$\deg g_j(x,y)\leq (p-j-1)n-3.\tag{5}$$

{\it Regularity conditions over
a nonsingular point of the branch curve.}
Consider our form
$$\omega =(\sum_{j=0}^{p-1}z^jg_j(x,y))\frac{dx\wedge dy}{z^{p-1}}$$
over a nonsingular point $(a,b)$ of
one of the components, $B_{i_0}$, of the branch curve.
Let $r_j$ be the order of zero (or of the pole if $r_j<0$) of the function
$g_j$ along the curve $B_{i_0}$, i.e., $g_j=\overline g_j\cdot h_{i_0}^{r_j}$
with $\overline g_j$
having neither poles nor
zeros along $B_{i_0}$.
Since $(a,b)$ is a nonsingular
point of $B_{i_0}$, we can assume that $h_{i_0}(x,y)$
and some function $g(x,y)$ are local analytic cootdinates in some
neighborhood $U$ of $(a,b)$ (denote them by $u$ and $v$).
So, over $U$ the surface
$X_G$ (after analytic change of
variebles) is isomorphic to the normalization $X_{G,loc}$ of the surface in
$\C ^3$ given by
$$z^{p}=u^{k_{i_0}}.$$
There is an analytic function $w$ in $X_{G,loc}$ such that $u=w^{p}$ and
$z=w^{k_{i_0}}$, and such that $w$ and $y$ are analytic coordinates in
$X_{G,loc}$. The differencial 2-form
$\omega$
considered above has the following form in the new coordinates
$$\omega =(\sum_{j=0}^{p-1}w^{jk_{i_0}}\overline g_j(x,y)w^{pr_j})
\frac{pw^{p-1}dw\wedge dv}{w^{(p-1)k_{i_0}}}.$$
It is easy to see that
$$j_1k_{i_0}+pr_{j_1}+p-1-(p-1)k_{i_0}\neq
j_2k_{i_0}+pr_{j_2}+p-1-(p-1)k_{i_0}$$
if $0<k_{i_0}<p$, $0\leq j_1, j_2\leq p-1$, and $j_1\neq j_2$.
Therefore,
all the rational functions $g_j(x,y)$
from $(4)$ are regular functions over a
nonsingular point $(a,b)$ of $B_{i_0}$ iff
$$jk_{i_0}+pr_{j}+p-1-(p-1)k_{i_0}\geq 0$$
Moreover, if $\omega$ is a regular form over $B_{i_0}$ then $r_j$ must be
greater than 0, since
for $0<k_{i_0}<p$, $0\leq j\leq p-1$, and $r_j\leq -1$, we obtain that
$$jk_{i_0}+pr_{j}+p-1-(p-1)k_{i_0}< 0.$$
From this it follows that if $\omega$ is a regular form then
all the rational functions $g_j(x,y)$ are regular functions
everywere in $\C ^2$ outside codimension 2,
and thus
$g_j(x,y)$ should be polynomials in $x$ and $y$.
Moreover, the polinomials $g_j(x,y)$ must be
divisible by $h^{r_j}_{i}(x,y)$,
where $r_j$ is the smallest integer satisfying the inequality
$$pr_j\geq (p-j-1)k_i -p+1. \tag{6} $$

{\it Regularity conditions over singular points of the branch curve.}
This is the only step where the singular points of $X_G$
are concerned. Let $\nu : \widetilde X_G\to X_G$
be the minimal resolution of singularities of $X_G$ and
$E$ be the exceptional divisor of $\nu$. Pick
a composition $\sigma : Y\to \Cp2$
of $\sigma$-processes with centers at
singular points of $B$ (and their preimages) such that
$\sigma ^{-1}\circ f\circ \nu (E_i)$ is a curve for each
irreducible component $E_i$ of $E$. Let $Z$ be the normalization of
$Y\times _{\Cp2} X_G$. Denote by $g: \widetilde X_G\to Z$ the
birational map induced by $\nu$ and $\sigma$. It follows from the
above choice of $\sigma$ that for any
$\omega \in H^0(Z\setminus \text{Sing}Z,
\Omega ^2_{Z\setminus \text{Sing}Z})$
its pull-back $g^*(\omega )$ is regular
at generic points of $E_i$ and, thus,
extends to a regular form on the whole
$\widetilde X_G$. Hence, $H^0(\widetilde X_G,
\Omega ^2_{\widetilde X_G})$
is isomorphic to
$H^0(Z\setminus \text{Sing}Z, \Omega ^2_{Z\setminus \text{Sing}Z})$.

Therefore, it remains to consider a 2-form $\omega $ written as in $(4)$ and
to find a criterion of its regularity on $Z\setminus\text{Sing}Z$.
It can be done by performing, step by step, the $\sigma$-processes
chosen above. Let us accomplish only the first step, since
it is sufficient for the calculation in our particular example
which follows.

Represent, once more, $X_G$ as normalization of the surface given by
$$z^p=h(x,y).$$
Denote by $r$ the order of zero of $h(x,y)$ at the point
$(0,0)$, $r=sp+q$, $0\leq q <p$, and perform  a
$\sigma$-process with center at this point.
In a suitable chart,
this $\sigma$-process $\sigma :\C ^2_{(u,v)}\to \C ^2_{(x,y)}$
is given by $x=u, \, \, y=uv$.
The normalization $Z_1$ of
$X_G\times _{\C ^2_{(x,y)}} \C ^2_{(u,v)}$
is birational to
the normalization
of the surface given by
$$w^p=u^q\overline h(u,v),$$
where $w=z/u^s$ and $\overline h(u,v)=h(u,uv)/u^r$.
We have
$$\align
\omega  & =(\sum_{j=0}^{p-1}z^jg_j(x,y))\frac{dx\wedge dy}{z^{p-1}}= \\
 & =(\sum_{j=0}^{p-1}w^j\overline g_j(u,v)u^{sj+s_j+1-s(p-1)})
\frac{du\wedge dv}{w^{p-1}},
\endalign
$$
where $s_j$ is the order of zero of $g_j(x,y)$ at $(0,0)$.
Applying $(6)$, we get necessary conditions for
the regularity of the pull-back of $\omega$
at generic points of the exceptional divisor:
the order of zero $s_j$ of each $g_j(x,y)$ at
singular point of the branch locus $B$ of order $r$ is
the smallest integer satisfying the inequality
$$ps_j\geq (p-j-1)r -2p+1. \tag{7}$$

{\bf Principal calculation.}
Here, we consider the Galois $G=\Z /5\Z \times \Z /5\Z$
covering of $\Cp2$ constructed in Section 2.
In this special case, the minimal desingularization
$\tilde X_G$ of $X_G$ is the induced covering of
$\Cp2$ blown up in the twelve multiple points of the
configuration.

According to the above algorithm, we should examine one
by one the six subgroups $H_1=((0,1))$, $H_2=((1,0))$,
$H_3=((4,1))$,
$H_4=((3,1))$,
$H_5=((2,1))$,
and $H_6=((1,1))$
of $G$ such that $G_i=G/H_i=\Z /5\Z$
(here we denote by $(\alpha)$ the cyclic subgroup
of $G$ generated by $\alpha\in G$).
The cyclic Galois coverings
$X_{G_i}\to\Cp2$ can be represented
as normalizations of the surfaces in $\C ^3$ given respectively by
$$
\align
z^5 & =l_1l_2l_3l_4^3l_5^3l_9 ,\\
z^5 & =l_1l_3l_4^3l_6l_7l_8^2l_9 ,\\
z^5 & =l_1^2l_2l_3^2l_4l_5^3l_6l_7l_8^2l_9^2,\\
z^5 & =l_1l_2^2l_3l_4^3l_5l_6^4l_7^4l_8^3l_9, \\
z^5 & =l_1l_2^4l_3l_4^3l_5^2l_6^2l_7^2l_8^4l_9, \\
z^5 & =l_2^2l_5l_6^3l_7^3l_8, \\
\endalign
$$

To calculate the geometric
genus of each of $\tilde X_{G_i}$
we find explicitly all the regular $2$-forms,
which we write as in (4).
We use the criteria $(5)-(7)$ (for convenience, we reproduce
in
Tables 2 - 4 below,
using notations
involved in $(5)-(7)$, the exact values of these bounds
evaluated in the case of $\Z/5$-coverings).
For $G_1$ we get the forms
$$
cl_4l_5l_9z\frac{dx\wedge dy}{z^4}, c\in \C;
$$
for $G_2$ we get
$$
(P_2l_4^2l_8+P_1l_4l_8z+cl_4z^2)\frac{dx\wedge dy}{z^4}
$$
where $P_i$ are polynomials in $x,y$ of degree $i$,
$c\in\C$, and
$$
%\{
p_{349}, p_{789}, p_{168},
p_{147}
%\}\subset
\in \{P_2=0\},
\, p_{147}\in\{P_1=0\};
$$
for $G_3$ we get
$$
(P_1l_1l_3l_5^2l_6l_8l_9+Q_1l_1l_3l_4l_5l_8l_9z+
c_1l_3l_8z^2+c_2z^3)\frac{dx\wedge dy}{z^4}
$$
where $P_i, Q_i$ are polynomials of degree $i$,
$c_k\in\C$,  and
$$
%\{
p_{123}, p_{789}
%\}\subset
\in\{P_1=0\},
%\{
p_{456}
%\}\subset
\in
\{Q_1=0\};
$$
for $G_4$ we get
$$
(P_2l_2l_4^2l_6^3l_7^3l_8^2+Q_2l_2l_4l_6^2l_7^2l_8z+
P_1l_4l_6l_7l_8z^2+Q_1z^3)\frac{dx\wedge dy}{z^4}$$
with
$$
%\{
p_{159}, p_{123}
%\}\subset
\in
\{P_2=0\},
%\{
p_{348}
%\}\subset
\in
\{Q_2=0\},
%\{
p_{267}
%\}\subset
\in
\{P_1=0\},
%\{
p_{267}
%\}\subset
\in
\{P_1=0\};
$$
for $G_5$ we get
$$
(P_2l_2^3l_4^2l_5l_6l_7l_8^3+P_1l_2^2l_4l_5l_6l_7l_8^2z+
Q_2l_2l_4l_8z^2+Q_1z^3)\frac{dx\wedge dy}{z^4}
$$
with
$$
%\{
p_{159}, p_{369}, p_{357}
%\}\subset
\in
\{P_2=0\},
%\{
p_{357}, p_{267}, p_{258}
%\}\subset
\in
\{Q_2=0\},
%\{
p_{258}
%\}\subset
\in
\{Q_1=0\};
$$
and, finally, for $G_6$,
$$
cl_2l_6l_7z\frac{dx\wedge dy}{z^4}.
$$
Therefore,
$p_g(\tilde X_{G_1})=1$, $p_g(\tilde X_{G_2})=5$,
$p_g(\tilde X_{G_3})=5$, $p_g(\tilde X_{G_4})=13$,
$p_g(\tilde X_{G_5})=11$, and $p_g(\tilde X_{G_6})=1$.

%Format: plain
$$\phantom{aaa}$$

\centerline{\vbox{%
\def\\{\cr\noalign{\hrule}}
\def\plain#1{\omit\hss#1\hss\srule}
\def\margin{\kern12pt}                 %% Adjust cell margins here
\def\srule{\vrule height9pt depth4pt} %% Adjust row height and depth here
\halign{%
 \vrule\margin\hss#\hss\margin\srule&&    %% First entry centered
 \margin\hss$#$\margin\vrule\cr           %% Others flushed right
 \noalign{\hrule}
% Table starts here
 $\deg g_j\le$ &\plain{$n=1$}&\plain{$n=2$}&\plain{$n=3$}&\plain{$n=4$}\\ %% Column headers
 $j=0$ &  1 & 5 &  8 & 13\\
 $j=1$ & 0 & 3 & 6 & 9 \\
 $j=2$ & 0 & 1 & 3 & 5\\
 $j=3$ & 0 & 0 & 0 & 1\\
 $j=4$ & 0 & 0 & 0 & 0\\
% Table ends here
 \crcr}}}

 %Format: plain
$$
%\phantom{aaa}
\text{
%Fig.
Table 2}$$

\centerline{\vbox{%
\def\\{\cr\noalign{\hrule}}
\def\plain#1{\omit\hss#1\hss\srule}
\def\margin{\kern12pt}                 %% Adjust cell margins here
\def\srule{\vrule height9pt depth3pt} %% Adjust row height and depth here
\halign{%
 \vrule\margin\hss#\hss\margin\srule&&    %% First entry centered
 \margin\hss$#$\margin\vrule\cr           %% Others flushed right
 \noalign{\hrule}
% Table starts here
 $
 %\deg
 r_j\ge$ &\plain{$k_i=1$}&\plain{$k_i=2$}&\plain{$k_i=3$}&\plain{$k_i=4$}\\
 %% Column headers
 $j=0$ &  0 & 1 &  2 & 3\\
 $j=1$ & 0 & 1 & 1 & 2 \\
 $j=2$ & 0 & 0 & 1 & 1\\
 $j=3$ & 0 & 0 & 0 & 0\\
 $j=4$ & 0 & 0 & 0 & 0\\
% Table ends here
 \crcr}}}

 %Format: plain
$$\text{Table 3}$$

\centerline{\vbox{%
\def\\{\cr\noalign{\hrule}}
\def\plain#1{\omit\hss#1\hss\srule}
\def\margin{\kern12pt}                 %% Adjust cell margins here
\def\srule{\vrule height9pt depth3pt} %% Adjust row height and depth here
\halign{%
 \vrule\margin\hss#\hss\margin\srule&&    %% First entry centered
 \margin\hss$#$\margin\vrule\cr           %% Others flushed right
 \noalign{\hrule}
% Table starts here
 $
 %\deg
 s_j\ge$ &\plain{$r=2$}&\plain{$r=3$}&\plain{$r=4$}&\plain{$r=5$}
 &\plain{$r=6$}&\plain{$r=7$}&\plain{$r=8$}&\plain{$r=9$}&\plain{$r=10$}
 &\plain{$r=11$}&\plain{$r=12$}\\ %% Column headers
 $j=0$ &  0 & 1 &  2 & 3 & 3 & 4 & 5 & 6 & 7 & 7 & 8\\
 $j=1$ & 0 & 0 & 1 & 2 & 2 & 3 & 3 & 4 & 5 & 5 & 6\\
 $j=2$ & 0 & 0 & 0 & 1 & 1 & 1 & 2 & 2 & 3 & 3 & 3\\
 $j=3$ & 0 & 0 & 0 & 0 & 0 & 0 & 0 & 0 & 1 & 1 & 1\\
 $j=4$ & 0 & 0 & 0 & 0 & 0 & 0 & 0 & 0 & 0 & 0 & 0\\
% Table ends here
 \crcr}}}

$$\text{Table 4}$$

Now, from Proposition 5.1 it follows that
$p_g(\tilde X_{G})=36$.
On the other hand, by Noether formula, $(c_1^2+c_2)/12=37$ and, thus,
$$q(\tilde X_G)=
\dim H^0(\tilde X_G,\Omega ^1_{\tilde X_G})=(c_1^2+c_2)/12-1-p_g=0,$$
i.e., $\tilde X_G$ is a regular surface.

\end
\widestnumber\key{DIK5}
\Refs

\ref\key BPV
\by
W.~Barth, C.~Peters, A.~Van de Ven
\book Compact Complex Surfaces
\yr 1984
\publ Springer-Verlag
\endref


\ref \key Ca
\by F.~Catanese
\paper Moduli spaces of surfaces
and real structures
\jour preprint math. AG/0103071
\endref

\ref \key FM
\by R.~Friedman, J.W.~Morgan
\paper Algebraic surfaces and $4$-manifolds
\jour  Bull. A.M.S.
\vol   18
\yr    1988
%\issue 3
\pages 1--19
\endref

\ref \key F
\by A.~Fujiki
\paper Coarse Moduli Space for Polarized Compact K\"ahler Manifolds
\jour  Publ. RIMS (Kyoto Univ.)
\vol   20
\yr    1984
%\issue 3
\pages 977--1005
\endref

\ref \key G
\by H.~Grauert
\paper Ein Theorem der analytischen Garbentjeorie
und due Modulra\"ume komplexer Strukturen
\jour Publ. IHES
\vol 5
\yr 1960
\endref

\ref \key I
\by M.-N.~Ishida
\paper The Irregularities of Hirzebruch's examples
of Surfaces of General Type with $C_1^2=3C_2$
\jour Math. Ann.
\vol 262
\yr 1983
\pages 407--420
\endref

\ref\key H \by  F.~Hirzebruch
\paper Arrangements of lines and algebraic surfaces
\inbook Arithmetics and Geometry
\vol II, Prog. Math. 36
\pages 113-140 \publ Birkh\'auser \yr 1983
\endref


\ref\key KK
\by  V. Kharlamov and Vik.S. Kulikov
\paper On real structures of rigid surfaces
\jour preprint AG/0101098 , to appear in Izvestiya RAN: Mathematics
\endref

\ref\key KK2
\by  V. Kharlamov and Vik.S. Kulikov
\paper Diffeomorphisms, isotopies and braid monodromy
factorizations of cuspidal curves
\jour preprint AG/0104021, to appear in C.R.A.S
\endref

\ref \key K
\by Vik.S.~Kulikov
\paper On Chisini's Conjecture
\jour Izvestiya Math.
\vol 63
\yr 1999
\pages 1139--1170
\endref

\ref\key KN
\by
S.~Kobayashi, K.~Nomizu
\book Foundations of differential geometry
\yr 1963
\publ Interscience publishers
\endref



\ref \key Ko
\by D.~Kotschick
\paper Orientations and geometrisations of compact complex surfaces
\jour Bull. London Math. Soc.
\vol 29
\yr 1997
\issue 2
\pages 145--149
\endref

\ref \key LeB
\by C.~LeBrun
\paper Diffeomorphisms, symplectic forms and Kodaira fibrations
\jour  Geometry \& Topology
\vol   4
\yr    2000
%\issue 3
\pages 451--456
\endref

\ref \key LL
\by T.-J.~Li, A.-K.~Liu
\paper Uniqueness of symplectic canonical class, surface cone
and symplectic cone of $4$-manifolds with $b^+=1$
\jour preprint math. 
SG/ 0012048
\endref

\ref \key LS
\by D.~Lieberman, E.~Sernesi
\paper Semicontinuity of $L$-dimension.
\jour Math. Ann.
\vol 225
\issue 1
\pages 77--88
\endref

\ref\key Ma \by Manetti M.
\paper  On the moduli space of diffeomorphic algebraic surfaces
\jour Invent. Math. \vol 143 (1) \pages  29--76 \yr 2001
\endref

\ref \key MT
\by C.T.~McMullen, C.H.~Taubes
\paper $4$-manifolds with inequivalent symplectic forms and
$3$-manifolds with inequivalent fibrations
\jour  Math. Research Let.
\vol   6
\yr    1999
%\issue 3
\pages 681--696
\endref

\ref\key Mu \by  D.~Mumford
\paper An algebraic surface with $K$ ample,
$K^2=9$, $p_g=q=0$
\jour Amer. J. Math.
\vol 101
\pages 233--244
\yr 1979  \endref


\ref \key R
\by Y.B.~Ruan
\paper Symplectic topology on algebraic $3$-manifolds
\jour  J. Diff. Geom.
\vol   39
\yr    1994
%\issue 3
\pages 215--227
\endref

\ref \key Shi
\by G.~Shimura
\paper On the real points of an arithmetic quotient of a bounded
symmetric domain
\jour  Math. Ann.
\vol   215
\yr    1975
%\issue 3
\pages 135--164
\endref

\ref\key S
\by  Y.-T.~Siu
\paper The complex-analyticity of
harmonic maps and the strong rigidity of compact K\"ahler manifolds
\jour Ann. Math
\vol 112 \pages
73--111 \yr 1980
\endref

\ref \key Sm
\by I.~Smith
\paper On moduli spaces of symplectic forms
\jour preprint math. SG/0012096
\endref

\ref \key W
\by C.T.C.~Wall
\paper On simply-connected 4-manifolds
\jour J. London Math. Soc.
\vol 39
\yr 1964
\pages 141--149
\endref

\ref \key Z
\by A.~V.~Zubr
\paper Classification of simply connected six-dimensional spin manifolds
\jour Izv. Akad. Nauk SSSR Ser. Mat.
\vol 39
\issue 4
\yr 1975
\pages 839--859
\endref
